\documentclass{article}
\title{On the twin-width of near-regular graphs}
\author{Irene Heinrich\thanks{Technische Universität Darmstadt, Germany}\\[1mm] Simon Raßmann\footnotemark[1] \and Ferdinand Ihringer\thanks{Southern University of Science and Technology (SUSTech), Shenzhen, China}\\[1mm] Lena Volk\footnotemark[1]}
\usepackage[utf8]{inputenc}
\usepackage[T1]{fontenc}
\usepackage{lmodern, microtype}

\usepackage{hyperref}
\hypersetup{unicode}

\usepackage{enumitem}

\usepackage{bm}

\makeatletter
\g@addto@macro\bfseries{\boldmath}
\makeatother

\usepackage{amsmath, amssymb, amsfonts, mathtools, stmaryrd, amsthm, thm-restate}
\usepackage{cleveref}
\usepackage{cancel}
\usepackage{tikz, pgf}
\usetikzlibrary{calc, cd, shapes, fit}
\SetSymbolFont{stmry}{bold}{U}{stmry}{m}{n}

\theoremstyle{definition}
\newtheorem{definition}{Definition}[section]

\newtheorem{observation}[definition]{Observation}

\theoremstyle{plain}
\newtheorem{lemma}[definition]{Lemma}

\newtheorem{theorem}[definition]{Theorem}	
\newtheorem{corollary}[definition]{Corollary}
\newtheorem{conjecture}[definition]{Conjecture}
\newtheorem{claim}{Claim}[definition]

\newenvironment{claimproof}{\proof[\(\ulcorner\)\nopunct]}{\endproof}

\makeatletter
\let\@vareps\varepsilon \let\varepsilon\epsilon \let\epsilon\@vareps
\let\phi\varphi
\makeatother

\DeclareMathOperator{\srg}{srg}
\DeclareMathOperator{\girth}{girth}

\tikzstyle{vertex}=[circle,draw,inner sep=0mm, minimum size=1.5mm, fill=black]
\tikzstyle{empty}=[]
\tikzstyle{edge}=[draw,thick]
\tikzset{>=stealth}

\usepackage{nicefrac}
\usepackage{pict2e,picture}

\newcommand{\dotcup}{\mathbin{\dot\cup}}

\renewcommand{\subset}{\subseteq}

\newcommand{\N}{\mathbb{N}}
\newcommand{\Z}{\mathbb{Z}}

\newcommand{\feq}{\preceq}

\DeclareMathOperator{\rdeg}{red-deg}
\newcommand{\maxdeg}{\Delta}
\newcommand{\maxrdeg}{\Delta_{\mathrm{red}}}
\newcommand{\red}{\operatorname{red}}
\DeclareMathOperator{\ord}{ord}

\DeclareMathOperator{\tww}{tww}
\DeclareMathOperator{\stww}{stww}
\DeclareMathOperator{\lb}{lb}
\DeclareMathOperator{\Aut}{Aut}
\DeclareMathOperator{\Cay}{Cay}
\DeclareMathOperator{\ls}{LS}

\makeatletter
\DeclareRobustCommand{\Arrow}[1][]{%
	\check@mathfonts
	\if\relax\detokenize{#1}\relax
	\settowidth{\dimen@}{$\m@th\rightarrow$}%
	\else
	\setlength{\dimen@}{#1}%
	\fi
	\sbox\z@{\usefont{U}{lasy}{m}{n}\symbol{41}}%
	\begin{picture}(\dimen@,\ht\z@)
		\roundcap
		\put(\dimexpr\dimen@-.7\wd\z@,0){\usebox\z@}
		\put(0,\fontdimen22\textfont2){\line(1,0){\dimen@}}
	\end{picture}%
}
\makeatother

\begin{document}
	\maketitle
		\begin{abstract}
	Twin-width is a recently introduced graph parameter based on the repeated contraction of near-twins.
	It has shown remarkable utility in algorithmic and structural graph theory, as well as in finite model theory---particularly since first-order model checking is fixed-parameter tractable
	when a witness certifying small twin-width is provided.
	However, the behavior of twin-width in specific graph classes, particularly cubic graphs, remains poorly understood.
	While cubic graphs are known to have unbounded twin-width, no explicit cubic graph of twin-width greater than~\(4\) is known.
	
	This paper explores this phenomenon in regular and near-regular graph classes.
	We show that extremal graphs of bounded degree and high twin-width are asymmetric,
	partly explaining their elusiveness. Additionally, we establish bounds for circulant and $d$-degenerate graphs,
	and examine strongly regular graphs, which exhibit similar behavior to cubic graphs.
	Our results include determining the twin-width of Johnson graphs over 2-sets, and cyclic Latin square graphs.
\end{abstract}
		\section{Introduction}
Since its introduction in \cite{tww1}, twin-width has rapidly gained prominence in algorithmic and structural graph theory, and finite model theory.
Most significantly, classes of bounded twin-width admit fixed-parameter tractable first order model checking when given a witness that the graph has small twin-width.
Moreover, bounded twin-width precisely coincides with tractable first order model checking on many hereditary graph classes, including
ordered graphs \cite{tww4}, permutation graphs \cite{tww1}, and interval graphs \cite{tww8_conference}.
However, classes of bounded twin-width also possess many further nice structural and algorithmic properties like smallness \cite{tww2}, an \(O(\log n)\)-bit adjacency labeling scheme \cite{tww2}, and \(\chi\)-boundedness \cite{tww3_conference}.
Since so many classes of graphs have bounded twin-width---including, among others, all classes of bounded tree-width, bounded rank-width, all classes excluding a minor, unit interval graphs, \(K_t\)-free interval graphs, and map graphs \cite{tww1}---these results generalize and unify many known results about these smaller graph classes.

Despite its broad applicability to graph classes, understanding twin-width for specific graphs remains challenging. Upper bounds on twin-width can be established by constructing contraction sequences, but almost all lower bounds rely on demonstrating the absence of near-twins---essentially proving that already the first contraction introduces significant errors.
This approach, however, fails for the class of cubic graphs. While cubic graphs
are known to have unbounded twin-width through counting arguments \cite{tww2},
constructing explicit examples with high twin-width has proven surprisingly difficult.
A reason for this difficulty is that due to their sparsity, all pairs of vertices of a cubic graph are near-twins in the sense that they only
have a bounded number of non-common neighbors, which implies that every cubic graph admits a partial contraction sequence of linear length and bounded
width (indeed, of width \(4\)).
\begin{figure}[b]
	\centering
	\begin{tikzpicture}[scale=0.47]
		\begin{scope}[xshift=0cm]
			\node[vertex] (0) at (90:2) {};
			\node[vertex] (2) at (90-72:2) {};
			\node[vertex] (6) at (90-2*72:2) {};
			\node[vertex] (4) at (90-3*72:2) {};
			\node[vertex] (1) at (90-4*72:2) {};
			
			\node[vertex] (3) at (90:1) {};
			\node[vertex] (9) at (90-1*72:1) {};
			\node[vertex] (7) at (90-2*72:1) {};
			\node[vertex] (8) at (90-3*72:1) {};
			\node[vertex] (5) at (90-4*72:1) {};
			
			\draw
			(0) -- (1) 
			(0) -- (2) 
			(0) -- (3) 
			(1) -- (4) 
			(1) -- (5) 
			(2) -- (6) 
			(2) -- (9) 
			(3) -- (7) 
			(3) -- (8) 
			(4) -- (6) 
			(4) -- (8) 
			(5) -- (7) 
			(5) -- (9) 
			(6) -- (7) 
			(8) -- (9);
		\end{scope}
		
		\begin{scope}[xshift=5cm]
			\coordinate (x) at (90:2);
			\node[vertex] (0) at ($(x)+(10:0.5)$) {};
			\node[vertex] (2) at ($(x)+(170:0.5)$) {};
			\node[vertex] (3) at ($(x)+(270:0.3)$) {};
			
			\node[vertex] (1) at (90-72:2) {};
			\node[vertex] (8) at (90-2*72:2) {};
			\node[vertex] (6) at (90-3*72:2) {};
			\node[vertex] (11) at (90-4*72:2) {};

			\node[vertex] (10) at (90:1) {};
			\node[vertex] (9) at (90-1*72:1) {};
			\node[vertex] (4) at (90-2*72:1) {};
			\node[vertex] (7) at (90-3*72:1) {};
			\node[vertex] (5) at (90-4*72:1) {};
			
			\draw
			(0) -- (1) 
			(0) -- (2) 
			(0) -- (3) 
			(1) -- (8) 
			(1) -- (9) 
			(2) -- (3) 
			(2) -- (11) 
			(3) -- (10) 
			(4) -- (5) 
			(4) -- (8) 
			(4) -- (10) 
			(5) -- (9) 
			(5) -- (11) 
			(6) -- (7) 
			(6) -- (8) 
			(6) -- (11) 
			(7) -- (9) 
			(7) -- (10) 
			;
		\end{scope}
		
		\begin{scope}[xshift=10cm]
			\node[vertex] (0) at (90:2) {};
			\node[vertex] (1) at (60:2) {};
			\node[vertex] (6) at (30:2) {};
			\node[vertex] (4) at (0:2) {};
			\node[vertex] (11) at (330:2) {};
			\node[vertex] (8) at (300:2) {};%
			\node[vertex] (3) at (270:2) {};%
			\node[vertex] (9) at (240:2) {};%
			\node[vertex] (10) at (210:2) {};%
			\node[vertex] (5) at (180:2) {};%
			\node[vertex] (7) at (150:2) {};%
			\node[vertex] (2) at (120:2) {};%
			
			\draw
			(0) -- (1) 
			(0) -- (2) 
			(0) -- (3) 
			(1) -- (6) 
			(1) -- (10) 
			(2) -- (7) 
			(2) -- (11) 
			(3) -- (8) 
			(3) -- (9) 
			(4) -- (6) 
			(4) -- (9) 
			(4) -- (11) 
			(5) -- (7) 
			(5) -- (8) 
			(5) -- (10) 
			(6) -- (7) 
			(8) -- (11) 
			(9) -- (10);
			
		\end{scope}
		
		\begin{scope}[xshift=15cm]
			
			\node[vertex] (1) at (90:2) {};
			\node[vertex] (0) at (30:2) {};
			\node[vertex] (2) at (330:2) {};
			\node[vertex] (9) at (270:2) {};
			\node[vertex] (5) at (210:2) {};
			\node[vertex] (11) at (150:2) {};
			
			\node[vertex] (6) at (90:1) {};
			\node[vertex] (3) at (30:1) {};
			\node[vertex] (7) at (330:1) {};
			\node[vertex] (10) at (270:1) {};
			\node[vertex] (8) at (210:1) {};
			\node[vertex] (4) at (150:1) {};
			
			\path[draw]
			(0) -- (1) 
			(0) -- (2) 
			(0) -- (3) 
			(1) -- (6) 
			(1) -- (11) 
			(2) -- (7) 
			(2) -- (9) 
			(3) -- (8) 
			(3) -- (10) 
			(4) -- (7) 
			(4) -- (10) 
			(4) -- (11) 
			(5) -- (8) 
			(5) -- (9) 
			(5) -- (11) 
			(6) -- (7) 
			(6) -- (8) 
			(9) -- (10) 
			;
		\end{scope}
		
		\begin{scope}[xshift=20cm]
			\node[vertex] (6) at (90:2) {};
			\node[vertex] (10) at (90-1*36:2) {};
			\node[vertex] (2) at (90-2*36:2) {};
			\node[vertex] (0) at (90-3*36:2) {};
			\node[vertex] (3) at (90-4*36:2) {};
			\node[vertex] (7) at (90-5*36:2) {};
			\node[vertex] (9) at (90-6*36:2) {};
			\node[vertex] (1) at (90-7*36:2) {};
			\node[vertex] (11) at (90-8*36:2) {};
			\node[vertex] (4) at (90-9*36:2) {};
			
			\node[vertex] (5) at ({-cos(54)}, 0) {};
			\node[vertex] (8) at ({cos(54)}, 0) {};
			
			\path[draw]
			(0) -- (1) 
			(0) -- (2) 
			(0) -- (3) 
			(1) -- (9) 
			(1) -- (11) 
			(2) -- (10) 
			(2) -- (11) 
			(3) -- (7) 
			(3) -- (8) 
			(4) -- (5) 
			(4) -- (6) 
			(4) -- (11) 
			(5) -- (8) 
			(5) -- (9) 
			(6) -- (7) 
			(6) -- (10) 
			(7) -- (9) 
			(8) -- (10) 
			;
		\end{scope}
	\end{tikzpicture}
	\caption{The five smallest cubic graphs of twin-width \(4\): The Peterson graph, the Tietze graph, the Triplex graph, the Twinplex graph and the Window graph.}
	\label{fig: small-tww4-cubic graphs}
\end{figure}

This hinders lower-bound arguments based on a shallow exploration of possible first contraction steps.
And in fact, no explicit construction of a cubic graph with twin-width greater than \(4\) is currently known (see Figure~\ref{fig: small-tww4-cubic graphs} for drawings of the five smallest cubic graphs of twin-width~4).
This stands in stark contrast to other graph classes where explicit constructions of graphs with unbounded twin-width,
such as the Rook's graphs~\cite{tww_products} or Paley graphs~\cite{bounds_on_tww}, are well-understood.

\paragraph*{Our results}
This paper investigates this phenomenon, exploring twin-width in regular and near-regular graph classes.

First, we show that extremal graphs of bounded degree and high twin-width are asymmetric. This result partly explains the difficulty in finding such graphs and sharply contrasts with the dense case, where
extremal examples often exhibit high symmetry (such as Rook's and Paley graphs, see~\cite{tww_products} and~\cite{bounds_on_tww}, respectively). Further, we give a quantitative strengthening of this result
by bounding the order of the automorphism group of near-extremal graphs of bounded degree and high twin-width. As part of this proof,
we show that the twin-width of circulant graphs is linearly bounded in their degree.

Second, we bound the twin-width of small subcubic graphs through an exhaustive computer search, showing that all subcubic graph on at most \(20\) vertices,
all cubic graphs on at most \(24\) vertices and all cubic graphs with girth at least \(6\) on at most \(28\) vertices have twin-width at most \(4\).

Third, we prove that every \(d\)-degenerate graph of order \(n\) has twin-width at most \(\sqrt{2dn}+2d\). While an \(O(\sqrt{dn})\)-bound
already follows from a more general bound in \cite{bounds_on_tww}, restricting to degenerate graphs makes our proof much more transparent
and deterministic.

Fourth, we explore strongly regular graphs, which exhibit a behavior similar to cubic graphs with respect to twin-width. Indeed, the minimum width
of a first contraction in a strongly regular graph is determined by its parameters. Though, most strongly regular graphs are pseudo-random \cite{pyber_pseudorandom}.
In this context, we determine the twin-width for several important graph families, including Johnson and Kneser graphs over 2-sets, cyclic Latin square graphs, and self-complementary edge- and vertex-transitive graphs (which, in particular, contain the class of all Peisert graphs). For all of these graphs, the twin-width agrees with the width of the first contraction.

\paragraph*{Related work}
The twin-width of cubic graphs was first proven to be unbounded in \cite{tww2} by showing that every class of bounded twin-width
is small, that is, contains only exponentially many isomorphism classes of graphs of order~\(n\), while the class of cubic graphs
is not small. The bound given there was significantly improved in \cite{tww_sparse_random_graphs}, where it is shown that cubic graphs
of order~\(n\) have twin-width at most \(n^{1/4+o(1)}\), and that this bound is attained with high probability by random cubic graphs.
In \cite{tww7_arXiv}, the fact that cubic graphs have unbounded twin-width was used to construct a finitely generated infinite group
such that the class of finite induced subgraphs of its Cayley graph has unbounded twin-width. However, no explicit construction
of such a group is known.

In the dense regime, general bounds on the twin-width of graphs were studied in \cite{bounds_on_tww}.
Among other results, the authors prove that every graph with~\(m\) edges has twin-width at most \(\sqrt{3m}+o(\sqrt{m})\),
which immediately yields that every \(d\)-degenerate graph \(G\) satisfies \(\tww(G)\leq \sqrt{3dn}+o(\sqrt{dn})\).
The proof of the bound in \cite{bounds_on_tww} is probabilistic: They show that for every \(m\)-edge graph, the vertex set can be
first contracted into at most \(\Theta(\sqrt{m})\) parts plus some exceptional vertices which can be handled using a randomized
contraction sequence. Then, only few parts remain which can be contracted arbitrarily.
We use  the more restricted setting of $d$-degeneracy to improve the constant factor in the bound,
and to give a much simpler and deterministic contraction sequence.

In the same paper, the authors also prove that Paley graphs of order~$n$ have twin-width \(\nicefrac{(n-1)}{2}\).
Further, the twin-width of Rook's graphs (and their generalization Hamming graphs) was determined in \cite{tww_products}.
This gives two families of strongly regular graphs, in each of which the twin-width was determined to be
the lower bound introduced by the first contraction step.

\paragraph*{Organization of the paper}
In \Cref{sec:preliminaries} we provide the necessary preliminaries and notation, including the formal definition of twin-width.
The focus of \Cref{sec:bounded_degree} are sparse graphs, including our results on extremal graphs of bounded degree, circulant graphs, small cubic graphs and degenerate graphs.
In \Cref{sec:strongly_regular} we explore strongly regular graphs and related families, presenting our findings on Johnson, Kneser, and Latin square graphs.

		\section{Preliminaries}\label{sec:preliminaries}
For $n \in \mathbb{N}$ we set $[n]\coloneqq \{1,2, \dots, n\}$.
We denote the set of $k$-element subsets of a set $V$ by $\binom{V}{k}$. For a prime power \(q\), we denote the finite field on~\(q\) elements by~\(\mathbb{F}_q\).

\paragraph*{Graph basics}
All graphs in this paper are finite and simple.
Let $G$ be a graph. We write $V(G)$ for the vertex set and $E(G)$ for the edge set of $G$, respectively.
The \emph{order} of $G$ is $|V(G)|$.
For $v \in V(G)$ we denote the \emph{degree} of $v$ in $G$ by $\deg_G(v)$ and we set \(\maxdeg(G) \coloneqq \max_{v\in V(G)}\deg_G(v)\).
We denote the \emph{(open) neighborhood} of~$v$ by~$N_G(v)$ and the \emph{closed neighborhood} of $v$ in $G$ by~$N_G[v]$.
For a set of vertices~\(S\subseteq V(G)\), we write \(N_G(S)\coloneqq\bigcup_{v\in S} N(S)\setminus S\) and \(N_G[S]\coloneqq \bigcup_{v\in S} N_G[v]\).
The \emph{girth} of~\(G\) is the length of the shortest cycle in~\(G\) and is denoted by~\(\girth(G)\).

\paragraph*{Partitions}
Let $S$ be a set.
A \emph{partition} of $S$ is a subset $\mathcal{P}$ of the power set of~$S$ with elements called ($\mathcal{P}$)-\emph{parts} such that $\bigcup_{P \in \mathcal{P}} P = S$ and each two distinct sets in $\mathcal{P}$ are disjoint.
If $\mathcal{P}$ and $\mathcal{Q}$ are partitions of the same set $S$, then $\mathcal{P}$ is a \emph{refinement} of $\mathcal{Q}$ (short: $\mathcal{P} \feq \mathcal{Q}$ ) if for each $P \in \mathcal{P}$ there is a $Q \in \mathcal{Q}$ such that~$P \subseteq Q$.
The \emph{discrete partition} of $S$ is $\{\{s\}\colon s \in S\}$ and the \emph{trivial partition} of $S$ is $\{S \}$.

\paragraph*{Graph isomorphisms}
An \emph{isomorphism} between two graphs $G$ and $H$ is a bijection $\varphi\colon V(G) \to V(H)$ which preserves adjacency and non-adjacency. An isomorphism from $G$ onto itself is an \emph{automorphism}.
A graph \(G\) is \emph{circulant} if there exists an automorphism \(\varphi\in\Aut(G)\) such that for each two vertices~$v$ and~$w$ of $G$ there exists some \(j\in\N\)
such that \(\varphi^j(v)=w\).
If~$G$ is isomorphic to its complement, then~$G$ is \emph{self-complementary}.

\paragraph{(Strongly) regular graphs}
We refer to~\cite{srgBrouwer}, \cite{cameron2004strongly}, and Chapter~10 of~\cite{DBLP:books/daglib/0037866} for literature on strongly regular graphs.
A graph \(G\) is \(d\)-regular if all vertices of~\(G\) have degree~\(d\).
If additionally $G$ is of order $n$  and there exist parameters~$\lambda$ and~$\mu$ such that for each two distinct vertices $u$ and $v$ of $G$ it holds that
\begin{equation*}
	\lvert N_G(u) \cap N_G(v)  \rvert = \begin{cases}
											\lambda &\text{if $u$ and $v$ are adjacent,}\\
											\mu &\text{otherwise},
										\end{cases}
\end{equation*}
then $G$ is a \emph{strongly regular graph} with parameter set $(n,d,\lambda, \mu)$ (short: $G$ is an $\srg(n,d,\lambda, \mu)$).
Every $\srg(n,d, \lambda, \mu)$ satisfies the following equality
\begin{equation}\label{eq: srg}
	(n-d-1)\mu = d(d-\lambda -1).
\end{equation}
 A \(3\)-regular graph is called \emph{cubic}, and
a graph of maximum degree at most \(3\) is \emph{subcubic}.

A strongly regular graph with parameters $(n, \frac{n-1}{2}, \frac{n-5}{4}, \frac{n-1}{4})$ for some $n \in \mathbb{N}$ is called a \emph{conference graph}. 
The \emph{Paley graph} \(\text{Paley}(q)\), for \(q \equiv 1 \pmod{4}\) a prime power, is the graph with vertex set \(\mathbb{F}_q\), where two vertices~$u$ and $v$ are adjacent precisely if~\(u-v\) is a nonzero square in \(\mathbb{F}_q\). Every Paley graph is a conference graph, see~\cite{srgBrouwer}. The \emph{Rook's graph} \(\text{Rook}(n)\) is the graph whose vertices are the cells of an \(n \times n\) chessboard, with edges joining two distinct cells precisely if they are either in the same row or the same column. 

\paragraph*{Degenerate graphs} 
A graph \(G\) is \emph{\(d\)-degenerate} if there exists a linear order~$v_1, v_2, \dots, v_n$ of~$V(G)$ such that
every vertex \(v_i\) has at most \(d\) left-neighbors, that is, neighbors among \(\{v_{1},\dots,v_{i-1}\}\).
This order is called a \emph{(\(d\)-)elimination order} of \(G\).\
\begin{observation}\label{obs:edges-of-degenerate-graphs}
	A \(d\)-degenerate graph of order \(n\) has at most \(d(n-1)\) edges.
\end{observation}
 
\paragraph*{Twin-width}
A \emph{trigraph} $G$ is a graph with edges colored either red or black.
Graphs are interpreted as trigraphs by coloring each edge black.
For $v \in V(G)$ the \emph{red degree} $\rdeg_G(v)$ is the number of red edges incident to~$v$.
The \emph{maximum red degree} of $G$ is $\maxrdeg(G) \coloneqq \max_{v \in V(G)}\rdeg_G(v)$.

Two disjoint vertex subsets $U$ and $W$ of $G$ are \emph{fully connected} if every 2-set \(\{u,w\}\) with \(u\in U\) and \(w\in W\) is a black edge. If no such pair is an edge, then \(U\) and \(W\) are \emph{disconnected}.
Given a partition~\(\mathcal{P}\) of \(V(G)\), we define the \emph{quotient graph} \(G/\mathcal{P}\) to be the trigraph with $V(G/\mathcal{P}) = \mathcal{P}$ such that two parts $U$ and $W$ of $\mathcal{P}$ are
\begin{align*}
	&\text{joined via a black edge} &\text{ if $U$ and $W$ are fully connected,}\\
	&\text{not adjacent} &\text{ if $U$ and $W$ are disconnected, and}\\
	&\text{joined via a red edge} &\text{ otherwise.}
\end{align*}
A \emph{partition sequence} of an order-$n$ trigraph \(G\) is a sequence \(\mathcal{P}_n, \mathcal{P}_{n-1}, \dots,\mathcal{P}_1\) of partitions 
of \(V(G)\), where \(\mathcal{P}_n\) is the discrete partition and for each \(i \in [n-1]\) the partition \(\mathcal{P}_i\) is obtained
by replacing two distinct parts
$P$ and $Q$ of \(\mathcal{P}_{i+1}\) by~\(P\cup Q\) (we call this a \emph{merge} or a \emph{contraction}). Equivalently, a partition sequence is given by the sequence
of trigraphs \(G/\mathcal{P}_i\)  called \emph{contraction sequence} of~$G$.
The \emph{width} of a contraction sequence (or its associated partition sequence) is the maximum red degree over all trigraphs~\(G/\mathcal{P}_i\).
A \emph{\(k\)-contraction sequence} is a contraction sequence of width at most~$k$.
The \emph{twin-width} \(\tww(G)\) of \(G\) is the minimal width of a contraction sequence of $G$.
If $G/{\mathcal{P}_n}, G/{\mathcal{P}_{n-1}}, \dots, G/{\mathcal{P}_1}$ is a contraction sequence of $G$,
then $G/{\mathcal{P}_n}, G/{\mathcal{P}_{n-1}}, \dots, G/{\mathcal{P}_i}$ is a \emph{partial contraction sequence} of $G$
for each $i \in [n]$, whose width is defined as for (complete) contraction sequences.

A variant of twin-width is \emph{sparse twin-width}~\cite{tww7_arXiv, tww_sparse_random_graphs},
which is defined as $\stww(G)\coloneqq\tww(G_{\red})$,
where \(G_{\red}\) is the trigraph obtained from \(G\) by coloring all edges red.
For graphs of bounded degree, twin-width and sparse twin-width differ by at most a constant.

\paragraph*{Lower bound for twin-width} A natural lower bound for the twin-width of a graph $G$ with at least two vertices follows from considering the minimum red degree obtained by just one merge:
\begin{equation} \label{eq: lb-is-a-lower-bound}
\tww(G) \geq
\min_{\substack{  \mathcal{P}_n, \mathcal{P}_{n-1}~\text{is a partial}\\ \text{contraction sequence of $G$}} } \maxrdeg(G/\mathcal{P}_{n-1}) \eqqcolon \lb_1(G).
\end{equation}
For $|V(G)| = 1$ we set $\lb_1(G) \coloneqq 0$.
A graph $G$ is \emph{$\lb_1$-collapsible} if~$\tww(G)=\lb_1(G)$.

\paragraph*{Groups}
All groups in this paper are finite, and we refer to \cite{finite_group_theory} as a standard textbook on finite group theory.
The \emph{order} of a group is the number of group elements.
We write \(S_n\) for the symmetric group on \(n\) elements and~\(\Z_n\) for the cyclic group
on \(n\) elements. A group whose order is a power of a prime~\(p\) is called a \(p\)-group.
Given a trigraph \(G\), we write \(\Aut(G)\) for the group of automorphisms of the underlying simple graph of \(G\),
thus, we do not require automorphisms to preserve edge colors.
A \emph{permutation group} on~\(G\) is a group~\(\Gamma\subseteq\Aut(G)\) with its natural action on \(V(G)\).
We call~$G$ \emph{vertex-transitive} or \emph{edge-transitive}  if \(\Aut(G)\) acts transitively on~$V(G)$ or~$E(G)$, respectively.
For a vertex~\(v\in V(G)\) and~\(\gamma\in\Gamma\), we write \(v^\gamma\coloneqq \gamma(v)\)
and~\(v^\Gamma\coloneqq\{v^{\gamma'}\colon\gamma'\in\Gamma\}\) for the \emph{orbit} of \(v\)
with respect to \(\Gamma\). The set of all orbits forms a partition of~\(V(G)\).
Given a group element \(\gamma\in\Gamma\), the \emph{order of \(\gamma\)} is the minimal \(n\geq 1\) such that~\(\gamma^n=1\)
and is denoted by \(\ord(\gamma)\). It coincides with the order of the subgroup \(\langle\gamma\rangle\subseteq \Gamma\)
generated by \(\gamma\).

A subgroup \(\Delta\subseteq\Gamma\) is \emph{normal} if for every \(\gamma\in\Gamma\) we have \(\gamma\Delta=\Delta\gamma\). 
A \emph{composition series} of a group \(\Gamma\) is a sequence
\(\Gamma=\Gamma_k\supseteq\Gamma_{k-1}\supseteq\dots\supseteq\Gamma_0=1\)
of subgroups of \(\Gamma\) of maximal length such that for all \(i\in[k]\), \(\Gamma_{i-1}\)
is a normal subgroup of \(\Gamma_i\).
If all quotient groups \(\Gamma_i/\Gamma_{i-1}\) of this composition series
are cyclic groups of prime order, then~\(\Gamma\) is \emph{solvable}.
By the Jordan-Hölder theorem, this does not depend on the composition series we chose.
Every finite \(p\)-group is solvable.

\paragraph*{Cayley Graphs} Let \(\Gamma\) be a group and \(S\) be an inverse-closed subset of \(\Gamma\). The \textit{Cayley graph} \(\Cay(\Gamma,S)\) is the graph with vertex set \(V(\Cay(\Gamma,S))=\Gamma\) and edge set \(E(\Cay(\Gamma,S)) = \{\{u,v\} \colon u^{-1}v \in S\}\). For an abelian group \(\Gamma\), we have~\(\{u,v\} \in E(\Cay(\Gamma,S))\) if and only if \(\{u^{-1},v^{-1}\} \in E(\Cay(\Gamma,S))\).

		\section{Twin-width of sparse graphs}\label{sec:bounded_degree}
\subsection{Extremal graphs of high twin-width}
In this section, we study extremal graphs of bounded degree and high twin-width,
that is, bounded-degree graphs with a minimal number of vertices for their twin-width.
For graphs of bounded degree, the notions of twin-width and sparse twin-width are functionally equivalent~\cite{tww7_arXiv, tww_sparse_random_graphs} via the bounds
\begin{equation}\label{eq:tww-vs-stww}
\tww(G)\leq \stww(G)\leq \tww(G)+\maxdeg(G).
\end{equation}
Our first result is that if the sparse twin-width of \(G\) is large enough, the first inequality of~\eqref{eq:tww-vs-stww} becomes an equality.
\begin{lemma}\label{lem:tww=stww}
If $G$ is a trigraph,
then
\[\stww(G)>\maxdeg(G)^2~\text{implies}~\stww(G)=\tww(G).\]
\begin{proof}
Let \(\mathcal{P}_n,\dots,\mathcal{P}_1\) be a partition sequence for \(G\) of width \(\tww(G)\).
We show that also \(\maxdeg(G_{\red}/\mathcal{P}_i)\leq \tww(G)\), which proves that \(\stww(G)\leq\tww(G)\).
Indeed, fix some \(i\in[n]\) and a part~\(P\in\mathcal{P}_i\).
If no other part of \(\mathcal{P}_i\) is fully connected to \(P\) in \(G\), this implies~\({\rdeg_{G_{\red}/\mathcal{P}_i}(P)=\rdeg_{G/\mathcal{P}_i}(P)\leq\tww(G)}\). Thus, assume
\(P\) is fully connected to some other part \(Q\) of \(\mathcal{P}_i\),
which implies~\({|P|\leq\maxdeg(G)}\). Moreover, since every vertex in \(P\) must be incident to \(Q\),
the number of edges from \(P\) to \(V(G)\setminus Q\) is at most \(|P|(\maxdeg(G)-1)\).
Thus, we find
\begin{align*}
\rdeg_{G_{\red}/\mathcal{P}_i}(P)
&\leq 1+|P|(\maxdeg(G)-1)
\leq 1+\maxdeg(G)^2-\maxdeg(G)\\
&\leq \stww(G)-\maxdeg(G)
\leq\tww(G),
\end{align*}
where the last inequality follows from \eqref{eq:tww-vs-stww}.
\end{proof}
\end{lemma}

\begin{corollary}\label{cor:tww_independent_of_edge_coloring}
For every trigraph \(G\) with underlying simple graph \(H\), if \({\stww(G)>\maxdeg(G)^2}\), then
\(\tww(G)=\tww(H)\).
\begin{proof}
Since \(\stww(G)=\stww(H)\), the claim follows from \Cref{lem:tww=stww}.
\end{proof}
\end{corollary}

In the following, we show that extremal subcubic graphs of high twin-width have girth at least \(5\), which can be used to diminish a combinatorial explosion in computer searches for subcubic graphs of high twin-width.
\begin{lemma}\label{lem:extremal-subcubic-graphs-girth5}
For every \(k>9\), every subcubic graph \(G\) with \(\tww(G)\geq k\) of minimal order
has girth at least \(5\).
\begin{proof}
By \Cref{cor:tww_independent_of_edge_coloring}, it suffices to show that every subcubic
graph \(G\) of girth less than \(5\) admits a partial \(9\)-contraction sequence
to some smaller subcubic trigraph \(H\). This implies that \(\tww(H)\geq \tww(G)\),
as otherwise, we can contract~\(H\) to obtain a \(\max\{\tww(H),9\}\)-contraction sequence of~\(G\).

If \(G\) contains a triangle, we contract the three vertices of the triangle
arbitrarily to a single vertex. This creates red degree at most \(3\), and the resulting trigraph is still subcubic.

If \(G\) contains a \(4\)-cycle, we can similarly contract a matching in the \(4\)-cycle.
The resulting trigraph will again be subcubic, but the intermediate trigraph obtained
after contracting one matching edge might have a vertex of degree~\(4\).
Again, the resulting trigraph must have twin-width at least \(\tww(G)\), which means that~\(G\) was not minimal.
\end{proof}
\end{lemma}

Note that plausibly, extremal subcubic graphs of high twin-width have even larger girth:
Large girth implies that the subgraphs induced on small parts of a partition sequence are relatively sparse,
which forces more edges between distinct parts of the partition. This plausibly leads to higher maximum degree
of the quotient trigraph.

We now investigate how symmetries of bounded-degree graphs can be exploited
to obtain partial contraction sequences of small width.
Given a trigraph~\(G\) and a group \(\Gamma\subseteq\Aut(G)\)\footnote{Recall that by \(\Aut(G)\), we mean the automorphism group of the simple graph underlying the trigraph~\(G\).}, we write~\(G/\Gamma\) for the trigraph obtained from~\(G\) by contracting every orbit of \(\Gamma\) to a single vertex.
\begin{lemma}\label{lem:quotients_degree}
	Let \(G\) be a trigraph and \(\Gamma\subseteq\Aut(G)\). For every vertex \(v\in V(G)\), we have
	\(\deg_{G/\Gamma}(v^\Gamma)\leq \deg_G(v)\). In particular, \(\maxdeg(G/\Gamma)\leq\maxdeg(G)\).
\begin{proof}
	Two orbits \(v^\Gamma\) and \(w^\Gamma\) share an edge in \(G/\Gamma\) if and only if \(v\) has a neighbor in \(G\) which is contained in \(w^\Gamma\).
\end{proof}
\end{lemma}

Thus, it is interesting to consider for which groups \(\Gamma\) there exists a contraction sequence of small width
that contracts \(G\) to \(G/\Gamma\). Our first result is that this is possible for cyclic groups.
\begin{lemma}
	\label{lem:stww-contract-obits}
	If~\(G\) is a trigraph and~\(\varphi\in\Aut(G)\), then there exists a partial \(4\maxdeg(G)\)-contraction sequence
	which contracts \(G\) to \(G/\langle\varphi\rangle\).
	In particular, \(\tww(G) \leq \max \{4\maxdeg(G), \tww(G/\langle\varphi\rangle)\}\).
\end{lemma}
\begin{proof}
	Let \(\{v_1^{\langle\varphi\rangle}, v_2^{\langle\varphi\rangle}, \dots, v_\ell^{\langle\varphi\rangle}\}\)
	be the set of orbits of~\(V(G)\) under the action of~\(\langle\varphi\rangle\)
	and write \(o_i\coloneqq |v_i^{\langle\varphi\rangle}|\) for the order of these orbits.
	On each orbit~\(v_i^{\langle\varphi\rangle}\) for \(i \in [\ell]\) the automorphism~\(\varphi\) is a cyclic permutation of length~\(o_i\), i.e., every vertex in \(V(G)\) is of the form \(v_i^{\varphi^j}\) for a unique \(i\in[\ell]\)
	and \(j\in[o_i]\).
	We thus also write \(v_i^j\coloneqq v_i^{\varphi^j}\) and \(V(G)=\{v_i^j\colon i\in[\ell], j\in[o_i]\}\).
	With this notation~\(\varphi\) acts as~\(\varphi(v_i^j)=v_i^{j+1\bmod o_i}\).
	We call~\(P\subseteq V(G)\) \emph{consecutive} if \(P=\{v_i^{j\bmod o_i}\colon j\in J\}\)
	for some \(i\in[\ell]\), and some interval \(J\subseteq\Z\).
	
	Fix \(i\in[\ell]\), \(k\in\N\), and set \(m\coloneqq\left\lceil \frac{o_i}{k} \right\rceil\).
	We define a partition of~\(v_i^{\langle\varphi\rangle}\) into~\(m\) consecutive parts
	of length at most \(k\) by setting
	\[P_i^q\coloneqq\begin{cases}
		\{v_i^j\colon (q-1)k<j\leq qk\} &\text{if } q<m,\\
		\{v_i^j\colon (q-1)k<j\leq o_i\} &\text{if } q=m
	\end{cases}\]
	for all \(q\in[m]\).
	Denote by \(\mathcal{P}^k\) the partition of \(V(G)\) obtained by partitioning every orbit in this way.
	Note that all the partition classes~\(P_i^q\) for~\(i \in [\ell]\) and~\(q \in [m-1]\) have exactly~\(k\) elements,
	and the classes~\(P_i^{m-1}\) have exactly \(o_i \bmod k\) elements.
	Note further that all parts \(P_i^q\) are consecutive and contain the full orbit \(v_i^{\langle\varphi\rangle}\)
	if~\(k\geq o_i\).

	Since~\(\varphi\) is a cyclic permutation on all of the orbits,
	the neighborhood of a vertex~\(v_i^j\) in~\(G\) is given by
	\begin{equation}\label{eq:neighbourhood}
		N_G(v_i^j)=\varphi^j(N_G(v_i))=\left\{v_{i'}^{j+j'\bmod o_{i'}} \colon v_{i'}^{j'} \in N_G(v_i)\right\}.
	\end{equation}
	Thus, for a partition class~\(P_i^q\), we can write
	\[N_G(P_i^q)=\bigcup_{v_{i'}^{j'}\in N_G\left(v_i^{(q-1)k+1}\right)} \{v_{i'}^{j'+j\bmod o_{i'}}\colon j\in[|P_i^q|]\},\]
	which is a union of at most \(\maxdeg(G)\) consecutive sets of length at most \(k\).
	
	But every consecutive set of length at most \(k\) in some orbit \(v_i^{\langle\varphi\rangle}\) intersects
	at most three parts \(P_i^q\) since at most one part per orbit has length less than~\(k\).
	Thus, every part of \(\mathcal{P}^k\) has degree at most \(3\maxdeg(G)\) in \(G/\mathcal{P}^k\).
	
	Further, in order to contract \(G/\mathcal{P}^k\) to \(G/\mathcal{P}^{2k}\), we only contract pairs of consecutive parts.
	Thus, the neighborhood in \(G\) of every part in \(\mathcal{P}^{2k}\) is a union of at most \(\maxdeg(G)\) consecutive sets
	of length at most \(2k\). Such sets can each intersect at most four parts of \(\mathcal{P}^k\),
	which means that the total degree of the contracted part is still bounded by \(4\maxdeg(G)\) while contracting \(G/\mathcal{P}^{k}\) to~\(G/\mathcal{P}^{2k}\).
	
	Thus, we find a partial contraction sequence
	\[G=G/\mathcal{P}^1\to G/\mathcal{P}^2\to G/\mathcal{P}^4\to G/\mathcal{P}^8\to\dots\to G/\mathcal{P}^{2^{\lceil\log_2 \ord(\varphi)\rceil}}\]
	of width at most~\(4\maxdeg(G)\), which in total contracts \(G\) to \(G/\langle\varphi\rangle\).
\end{proof}

If \(G\) is a circulant graph, a more precise analysis of the above proof yields an even sharper bound:
\begin{lemma}\label{lem:circulant}
	If \(G\) is a circulant graph, then~\(\stww(G)\leq 3\maxdeg(G)+1\).
\begin{proof}
	Since \(G\) is circulant, there exists an automorphism \(\varphi\in\Aut(G)\) such that~\(G/\langle\varphi\rangle\) is the singleton graph.
	Thus, \Cref{lem:stww-contract-obits} yields a contraction sequence for \(G\) of width at most \(4\maxdeg(G)\). We argue
	that the contraction sequence constructed there actually has width at most \(3\maxdeg(G)+1\).
	Indeed, recall that we ordered \(V(G)\) using the automorphism \(\varphi\),
	and partitioned \(V(G)\) into intervals of length \(k\) and at most one interval \(P^m\) of length less than \(k\).
	Then, the neighborhood of one of these intervals is a union of at most \(\maxdeg(G)\) other intervals, each of which has length at most \(k\).
	Such an interval now either intersects at most two parts of the partition, or it intersects exactly three parts, one of which is \(P^m\).
	Thus, excluding \(P^m\), each interval intersects at most two parts, which means that the neighborhood in total intersects at most \(2\maxdeg(G)+1\)
	parts. If we contract two consecutive parts, the neighborhood of this new contracted part is a union of at most \(\maxdeg(G)\) intervals of length at most \(2k\).
	Using the same reasoning as before, this neighborhood intersects at most \(3\maxdeg(G)+1\) parts. This shows that the constructed contraction sequence
	has width at most \(3\maxdeg(G)+1\).
\end{proof}
\end{lemma}

We continue by using \Cref{lem:stww-contract-obits} to show that the twin-width-extremal graphs of bounded degree are asymmetric.
\begin{lemma}
	\label{lem:asym-quotient}
	Every trigraph \(G\) admits a partial \(4\maxdeg(G)\)-contraction sequence to some
	asymmetric trigraph \(H\) with \(\maxdeg(H)\leq\maxdeg(G)\).
\end{lemma}
\begin{proof}
	Let~\(H\) be a quotient of~\(G\) of maximum degree at most~\(\maxdeg(G)\) to which~\(G\) can be contracted while keeping the degree of all intermediate trigraphs bounded by \(4\maxdeg(G)\),
	and choose~\(H\) to be of minimal order (such a graph exists as~\(G\) itself is such a graph).
	We show that~\(H\) is asymmetric. Indeed, assume for a contradiction that this is not the case,
	and pick some non-trivial automorphism~\(\varphi\in\Aut(H)\).
	Applying \Cref{lem:stww-contract-obits} to the graph~\(H\), we can contract each of the \(\langle\varphi\rangle\)-orbits into a single
	vertex while keeping the degree bounded by~\(4\maxdeg(G)\).
	The resulting graph~\(H/\langle\varphi\rangle\) has maximum degree at most~\(\maxdeg(H)\leq\maxdeg(G)\) by \Cref{lem:quotients_degree}.
	This contradicts the minimality of~\(H\).
\end{proof}

\begin{corollary}\label{cor:extremal_trigraphs_asymmetric}
	Let~\(d\geq 3\) and~\(k>\max\{4d,d^2\}\). Then every graph of maximum degree at most~\(d\)
	and twin-width at least~\(k\) of minimal order is asymmetric.
\end{corollary}

\begin{proof}
	Let~\(G\) be a graph of minimal order with $\Delta(G) \leq d$ and $\tww(G) \geq k$.
	Assume that~\(G\) is not asymmetric. By Lemma~\ref{lem:asym-quotient} the graph~\(G\) admits a partial
	\(4d\)-contraction sequence to some strictly smaller trigraph \(H\) with maximum degree
	at most \(d\). 
	
	Since \(k\leq \tww(G)\leq\max\{4d,\tww(H)\}\), and~\(k>4d\), we have~\(\tww(H)\geq k\).
	Thus~\(\Delta(H) \leq d\) and~\(\tww(H) \geq k>d^2\).
	By \Cref{cor:tww_independent_of_edge_coloring}, the underlying simple graph of~\(H\) also has
	twin-width at least \(k\), which contradicts the minimality of~\(G\).
\end{proof}

Note that neither the bound on the girth in \Cref{lem:extremal-subcubic-graphs-girth5} nor asymmetry in \Cref{cor:extremal_trigraphs_asymmetric}
are true in general in the case of unbounded degree:
The extremal trigraphs in the dense and sparse settings
include red cliques and red stars, both of which have large automorphism groups.
When restricting to graphs, the minimal graphs
of twin-width 1, 2, 3 and 4 were determined in \cite{SAT_approach},
and all but one of them have non-trivial symmetries.
Maybe more strikingly, it is an open question in \cite{bounds_on_tww}
whether every graph has \(\tww(G)\leq\frac{|G|-1}{2}\),
a bound achieved by the Paley graphs. A positive resolution
would imply that Paley graphs are extremal for twin-width,
while also being edge-transitive and thus highly symmetric.

Next, we want to show that we can not only contract the orbits induced by a cyclic group,
but also those induced by arbitrary solvable subgroups of the automorphism group.
\begin{lemma}\label{lem:solvable_groups}
	Let \(G\) be a trigraph and \(\Gamma\subseteq\Aut(G)\) be solvable.
	Then there exists a partial \(4\maxdeg(G)\)-contraction sequence which contracts \(G\) to \(G/\Gamma\). 
	\begin{proof}
	Let \(\Gamma=\Gamma_k\supseteq\Gamma_{k-1}\supseteq \dots\supseteq \Gamma_0=1\) be a composition series of \(\Gamma\),
	i.e., a sequence of subgroups ending in the trivial group such that for all \(i\), \(\Gamma_i\subseteq\Gamma_{i+1}\) is normal and
	each quotient group \(\Gamma_{i+1}/\Gamma_i\) is a cyclic group of prime order.
	Using this composition series, the contraction of \(G\) to \(G/\Gamma\) can be split into a sequence of contractions
	\[G=G/\Gamma_0\to G/\Gamma_1\to G/\Gamma_2\to \dots\to G/\Gamma_k=G/\Gamma,\]
	where every intermediate trigraph again has maximum degree at most \(\maxdeg(G)\) by \Cref{lem:quotients_degree}.
	Thus, we only need to argue that for every \(i\in[k]\), the trigraph~\(G/\Gamma_{i-1}\)
	can be contracted to \(G/\Gamma_i\) while never exceeding degree~\(4\maxdeg(G)\).
	
	For this, consider the action of the group \(\Gamma_i/\Gamma_{i-1}\) on the trigraph \(G/\Gamma_{i-1}\) via~\mbox{\((v^{\Gamma_{i-1}})^{\gamma_i\Gamma_{i-1}}\coloneqq \left(v^{\gamma_i}\right)^{\Gamma_{i-1}}\)},
	which is well-defined since~\(\Gamma_{i-1}\) is normal in~\(\Gamma_i\).
	Further, only the cosets \(\gamma_i\Gamma_{i-1}\) with \(\gamma_i\in\Gamma_{i-1}\) act trivially,
	meaning~\(\Gamma_i/\Gamma_{i-1}\) is isomorphic to a subgroup of \(\Aut(G/\Gamma_{i-1})\).
	Since \(\Gamma_i/\Gamma_{i-1}\) is a cyclic group of prime order,
	this implies that we find some automorphism~\(\varphi\in\Aut(G/\Gamma_{i-1})\)
	such that \((G/\Gamma_{i-1})/\langle\varphi\rangle=G/\Gamma_i\).
	The orbits of this automorphism can be contracted using width at most \(4\maxdeg(G/\Gamma_{i-1})\) by \Cref{lem:stww-contract-obits},
	which is within our bound of~\(4\maxdeg(G)\) by \Cref{lem:quotients_degree}.
	\end{proof}
\end{lemma}
The above lemma implies that every vertex-transitive graph \(G\) with solvable automorphism group has twin-width at most \(4\maxdeg(G)\).

As an application of \Cref{lem:solvable_groups}, we prove a quantitative strengthening of \Cref{cor:extremal_trigraphs_asymmetric}.
Indeed, we show that for graphs with bounded degree, a large automorphism group implies that
the graph is in some sense far from being extremal for twin-width.
\begin{theorem}\label{thm: extremal trigraphs of bd degree}
	Let \(\ell\geq 0\) and let \(G\) be a graph of maximum degree~\(d\) and twin-width at least \(k>\max(4d,d^2)\).
	If every graph~\(H\) of maximum degree~\(d\) and order less than~\(|V(G)|-\ell\)
	has twin-width less than~\(k\), then \(\lvert\Aut(G)\rvert\leq O(\ell)^\ell\). 
\begin{proof}
	Let \(G\) be a graph as in the statement of the theorem.
	By \Cref{lem:solvable_groups}, we can find for every solvable subgroup \(\Gamma\) of \(\Aut(G)\)
	a partial \(4d\)-contraction sequence of \(G\) to \(H=G/\Gamma\).
	This implies \(\tww(H)\geq \tww(G)>d^2\), which by \Cref{cor:tww_independent_of_edge_coloring}
	and our assumption implies that \(|V(H)|\geq	 |V(G)|-\ell\).
	Hence, every solvable subgroup \(\Gamma\) of \(\Aut(G)\) must induce at least	\(|V(G)|-\ell\) orbits on~\(V(G)\).
	
	We proceed to show how the absence of a solvable subgroup with few orbits
	implies that all prime powers that divide the order of \(\Aut(G)\) must be small.
	\begin{claim}
		If \(p^e\) is a divisor of \(\lvert\Aut(G)\rvert\) for a prime \(p\) and some \(e\in\N\), then~\(e(p-1)\leq \ell\).
	\end{claim}
		
	\begin{claimproof}
		By the first Sylow theorem \cite{finite_group_theory}, we find a subgroup \(\Gamma\subseteq\Aut(G)\) of order~\(p^e\).
		Since every finite \(p\)-group is solvable \cite{finite_group_theory}, this implies
		that \(\Gamma\) induces at least~\(|V(G)|-\ell\) orbits on \(V(G)\).
		
		Let \(V(G)=O_1\dotcup\dots\dotcup O_t\) be the partition of \(V(G)\) into orbits with respect to \(\Gamma\).
		For every \(i\in[t]\), we get a restriction homomorphism \(\Gamma\to\Aut(G[O_i])\),
		whose image is again a \(p\)-group \(\Gamma_i\) of order \(p^{e_i}\).
		Thus, we get an embedding of \(p\)-groups
		\(\Gamma\to\prod_{i=1}^t\Gamma_i\), where the latter group induces the same orbit partition.
		By the orbit-stabilizer theorem, the order of every orbit divides the group order \(p^e\),
		hence for all \(i\in[t]\) we have \(|O_i|=p^{o_i}\) for some exponent \(o_i\in\N\).
		Thus, \(\Gamma_i\) embeds as a subgroup into \(S_{p^{o_i}}\), which implies that the order of~\(|\Gamma_i|\) divides \(p^{o_i}!\) by Lagrange's theorem \cite{finite_group_theory}.
		Since \(|\Gamma_i|=p^{e_i}\), the exponent \(e_i\) is bounded by the exponent of \(p\)
		in the prime decomposition of \(p^{o_i}!\).
		Since the product \(p^{o_i}!\) contains exactly \(p^{o_i}/p^j\) multiples
		of \(p^j\) for all \(j\leq o_i\), this yields
		\(e_i\leq p^{o_i-1}+p^{o_i-2}+\dots+p+1=\frac{p^{o_i}-1}{p-1}=\frac{|O_i|-1}{p-1}\).
		Rearranging this, we get \(e_i(p-1)\leq |O_i|-1\). Summing up these equations for all \(i\in[t]\) yields that~\(e(p-1)=\left(\sum_{i=1}^t e_i\right)(p-1)\leq |V(G)|-t\leq \ell\).
	\end{claimproof}
	
	The above claim restricts the prime powers that can appear in the prime factorization of \(\lvert\Aut(G)\rvert\).
	Indeed, it implies that the exponent of the prime \(p\) in this factorization is at most \(\frac{\ell}{p-1}\),
	and in particular, only primes up to \(\ell+1\) can appear at all.
	This gives an upper bound
	\[\lvert\Aut(G)\rvert\leq\prod_{p\leq \ell+1}p^{\frac{\ell}{p-1}}=\left(\prod_{p\leq \ell+1} p^{\frac{1}{p-1}}\right)^\ell,\]
	where the products range over all prime numbers up to \(\ell+1\).
	Taking the natural logarithm of both sides yields
	\[\ln(\lvert\Aut(G)\rvert)\leq \ell\sum_{p\leq \ell+1}\frac{\ln(p)}{p-1}=\ell\sum_{p\leq \ell+1}\frac{\ln(p)}{p}+\ell\sum_{p\leq \ell+1}\frac{\ln(p)}{p(p-1)}.\]
	By Merten's first theorem \cite{mertens}, the first sum is bounded by \(\ln(\ell+1)+2\), while the second sum is convergent and thus bounded independently of \(\ell\).
	This yields \(\ln(\lvert\Aut(G)\rvert)\leq \ell(\ln(\ell+1)+O(1))\) and thus
	\(\lvert\Aut(G)\rvert\leq O(\ell)^\ell\).
\end{proof}
\end{theorem}

We hope that the results obtained so far can eventually help with the construction of bounded-degree graphs
with high twin-width by pointing towards constructions that are less symmetrical in nature and possess high girth.

\subsection{Twin-width of small graphs of bounded degree}
We end our study of graphs of bounded degree with a computational result on the twin-width of small (sub-)cubic graphs. 

\begin{lemma}
	\label{lem:comput-small-cubic}
	The following graph classes have twin-width at most~4:
	\begin{enumerate}
	\item subcubic graphs of order at most \(20\),
	\item cubic graphs of order at most \(24\),
	\item cubic graphs of order at most \(28\) and girth at least \(6\).
	\end{enumerate}
\end{lemma}
\begin{proof}
	The list of subcubic graphs was generated using geng from version 2.8090 of the nauty package \cite{HPnauty, nauty}.
	The list of all such cubic graphs is available at~\cite{house_of_graphs_cubic, house_of_graphs} and was independently
	confirmed by genreg \cite{genreg}, minibaum \cite{minibaum} and snarkhunter~\cite{snarkhunter}.
	The twin-width of all these graphs was bounded by computer calculations using the heuristic solver GUTHM~\cite{guthm, guthm_software}.
	Only on the few graphs where the heuristic solver did not return a contraction sequence of width at most~\(4\) within our time bounds
	of between \(0.01\,\mathrm{ms}\) and \(0.5\,\mathrm{ms}\) depending on the graph size, we used version 0.0.3-SNAPSHOT
	of the exact twin-width solver hydraprime~\cite{hydraprime} to determine the exact twin-width.
\end{proof}

\subsection{Degenerate graphs}
In \cite{bounds_on_tww}, the authors give a probabilistic proof that every graph with at most \(m\) edges has twin-width at most \(\sqrt{3m}+o(\sqrt{m})\).
Combining this with Observation~\ref{obs:edges-of-degenerate-graphs} we obtain that every \(d\)-degenerate graph \(G\) satisfies~\({\tww(G)\leq \sqrt{3dn}+o(\sqrt{dn})}\).
For \(d\)-degenerate graphs we provide a simpler proof of a sharper twin-width bound in the following.

\begin{theorem}\label{thm:degeneracy}
Every \(d\)-degenerate graph \(G\) satisfies \(\tww(G)\leq\sqrt{2d|G|}+2d\).
\begin{proof}
We may assume that $G$ is $d$-degenerate but not $(d-1)$-degenerate with~\(V(G) = [n]\) and $d$-elimination order $1,2, \dots, n$.

We set \(k\coloneqq\left\lceil\sqrt{\frac{2n}{d}}\right\rceil\), which implies that \(d\sum_{i=1}^k i\geq n\).
Now, we partition~\([n]\) as follows:
For every \(i\in[dk]\), let \(I_i\subseteq[n]\) be the set containing
the \(\lceil\frac{i}{d}\rceil\) smallest elements in \([n]\setminus\bigcup_{j<i}I_i\),
or all elements in \([n]\setminus\bigcup_{j<i}I_i\) if there are at most \(d-1\) left.
Finally, we let \(\ell\in[dk]\) be the largest index for which \(I_\ell\) is non-empty. We obtain a partition $\{I_1, I_2, \dots, I_{\ell}\}$ of $[n]$.

We start by contracting each of the intervals into a single vertex as follows:
In each step, we choose the right-most interval \(I_i\) which is not yet contracted to a single part,
and in this interval contract the two rightmost parts.

Finally, we end up with the partition into intervals, which we contract arbitrarily.

Let \(\mathcal{P}\) be a partition along this partition sequence before two parts from distinct intervals are contracted. Then there is a vertex \(v\) such that all vertices~\(w\leq v\) lie in singleton parts,
while every two vertices \(w,w'>v\) lie in the same part if and only if they lie in the same interval \(I_i\). In particular, each part of~\(\mathcal{P}\) is convex with respect to the elimination ordering.

We first determine the \emph{red left-degree} of each part \(P\in\mathcal{P}\), that is,
the number of red edges joining \(P\) in the trigraph \(G/\mathcal{P}\) to parts \(Q\in\mathcal{P}\) 
that are smaller than \(P\) with respect to the elimination ordering.
If \(P\) is a singleton part, it has red left-degree \(0\). If \(P\) is not a singleton part,
the left-degree of every vertex in \(P\) is bounded by \(d\), which implies that the total
left-degree of~\(P\) is bounded by \(|P|d\). If \(P\subseteq I_i\), this is bounded by \(|I_i|d\).
Since every red edge in \(G/\mathcal{P}\) is incident to at least one non-trivial part and every
interval \(I_j\) contains at most one non-trivial part, the red right-degree of \(P\) is bounded by the number of intervals to the right of the interval \(I_i\) containing \(P\), including~\(I_i\) itself.
This number is bounded by \(d(k+1-|I_i|)\). Thus, the total red degree of \(P\) is bounded by
\(|I_i|d+d(k+1-|I_i|)=d(k+1)\).

After all intervals are contracted, there are only \(\ell\leq dk\) parts left. Thus, the red degree from this point on is also bounded by \(dk<d(k+1)\).
In total, this sequence has red degree at most
\(d(k+1)=d\left\lceil\sqrt{\frac{2n}{d}}\right\rceil+d\leq \sqrt{2dn}+2d\).
\end{proof}
\end{theorem}

When \(d\in\omega(\log n)\), then the bound in \Cref{thm:degeneracy} is asymptotically tight even for random graphs, which follows from \cite{tww_random_graphs}
together with the observation that random graphs in this regime asymptotically almost surely have maximum degree~\((1+o(1))d\)~\cite{bollobas_maxdegree}. 
However, when \(d\) is small or even fixed, we do not know whether it is tight.
If in addition to \(d\)-degeneracy, we require that \(\Delta(G)\leq2^{\sqrt[6]{\log(n)}}\), it follows from \cite[Lemma 2.8]{tww_sparse_random_graphs} that
\[\tww(G)\leq \Delta(G)\cdot n^{\frac{d-1}{2d-1}+o(1)},\]
which for fixed \(d\) is asymptotically smaller than the bound \(\sqrt{2dn}+2d\) we obtained.

		\section{Strongly regular graphs}\label{sec:strongly_regular}
Strongly regular graphs are a natural class of graphs to consider for twin-width for several reasons.
First, the red degree created by a single contraction is fully determined by the parameters of the strongly regular graph
and only depends on whether the contracted vertices are adjacent or non-adjacent~\cite{SAT_approach}.
Taking the minimum of these two options, we get
\[\lb_1(\srg(n,d,\lambda,\mu))=\min(2(d-\lambda-1),2(d-\mu)).\]
This natural lower bound for strongly regular graphs is often quite large compared to their order.
In particular, Paley graphs, which are strongly regular with parameters \(\left(n,\frac{n-1}{2},\frac{n-5}{4},\frac{n-1}{4}\right)\),
have twin-width \(\nicefrac{n-1}{2}\) and it is unknown whether \(n\)-vertex graphs with larger twin-width exist \cite{bounds_on_tww}.

A second reason for the interest in strongly regular graphs is that while the natural lower bound to their twin-width
is easy to understand, large strongly regular graphs behave like random graphs. For example,
apart from Latin square graphs and Steiner graphs, all strongly regular graphs
have a large eigenvalue gap, which means that they behave pseudo-randomly \cite{pyber_pseudorandom}.

In \cite{bounds_on_tww}, the authors showed that every \(n\)-vertex graphs has \(\lb_1(G)\leq\frac{n-1}{2}\).
We recall their proof and characterize those graphs for which equality occurs.

\begin{theorem}\label{thm:bounds_on_lb1}
If \(G\) is a graph of order \(n\), then \(\lb_1(G)\leq\frac{n-1}{2}\) with equality if and only if \(G\) is a conference graph.
\begin{proof}
Pick a contraction pair \(\{u,v\}  \in \binom{V}{2}\) uniformly at random. For every vertex~\(w\in V(G)\),
let \(R_w\) be the event that the contraction of \(\{u,v\}\) creates a red edge to~\(w\).
Since the contraction pairs that create a red edge to \(w\) are exactly those that contain one neighbor and one non-neighbor of \(w\) we obtain
\begin{equation}\label{eq:conference:red_degree_to_single_vertex}
\mathbb{P}(R_w)=\frac{\deg(w)\cdot(n-1-\deg(w))}{\binom{n}{2}}\leq\frac{\left(\frac{n-1}{2}\right)^2}{\binom{n}{2}}=\frac{n-1}{2n}.
\end{equation}

This yields
\begin{equation}\label{eq:conference:total_red_degree}
\mathbb{E}(\maxrdeg(G/\{u,v\}))=\sum_{w\in V(G)} \mathbb{P}(R_w)\leq n\cdot\frac{n-1}{2n}=\frac{n-1}{2}.
\end{equation}
In particular, there is a set $\{u,v\} \in \binom{V}{2}$ with \(\maxrdeg(G/\{u,v\})\leq\frac{n-1}{2}\),
which implies \(\lb_1(G)\leq\frac{n-1}{2}\).

If \(\lb_1(G)=\nicefrac{(n-1)}{2}\), then
the inequalities in (\ref{eq:conference:red_degree_to_single_vertex}) and (\ref{eq:conference:total_red_degree})
become equalities and \(G\) is \(\nicefrac{(n-1)}{2}\)-regular.
Further, the expected red degree after one random contraction is equal to the minimal red degree after any contraction and, hence,
every possible first contraction yields a red degree of exactly~\(\nicefrac{(n-1)}{2}\).
We obtain that $G$ is an $\srg(n,\frac{n-1}{2}, \frac{n-5}{4}, \frac{n-1}{4})$, which is the characterizing parameter set of conference graphs.
\end{proof}
\end{theorem}

Now we turn again to general strongly regular graphs and show that---at least early in the contraction sequence---contraction sequences of small width
are quite structured.
\begin{lemma}
	\label{lem:contrac-seq-is-struct}
Let \(G\) be an \(\srg(n,d,\lambda,\mu)\).
If \(\mathcal{P}_n,\mathcal{P}_{n-1},\dots,\mathcal{P}_{k}\) is a partial \(\lb_1\)-contraction sequence of \(G\), then \(|P|\leq 2\)
for all \(P\in\mathcal{P}_i\) with \(i\geq n-\left\lceil\frac{\lb_1(G)}{2}\right\rceil\).
\begin{proof}
	Let \(\mathcal{P}_n,\mathcal{P}_{n-1}, \dots, \mathcal{P}_k\) be a partial partition sequence of $G$ of width at most~\(\lb_1(G)\)
	and assume that \(\mathcal{P}_k\) is the first partition in this sequence that contains a part \(P\) of size greater than two.
	Further, denote by \(\mathcal{P}_P\) the partition of \(V(G)\) whose only non-singleton part is \(P\).
	We first give a bound for $\rdeg_{G/\mathcal{P}_P}(P)$.
	
	First assume that \(P=\{u,v,w\}\) contains exactly three vertices. We show that \(\rdeg_{G/\mathcal{P}_P}(P)\geq\frac{3}{2}\lb_1(G)-1\). The red neighbors of \(P\) in \(G/\mathcal{P}_P\) are the neighbors of either \(u\), \(v\) or \(w\), without
	the common neighbors of all three vertices and without the three vertices themselves. Thus, assume that the three vertices~$u$,~$v$, and~$w$ have exactly \(\nu\) common neighbors and further assume for simplicity that \(u\), \(v\) and \(w\) are pairwise non-adjacent.
	By the inclusion-exclusion principle, we have
	\(|N(\{u,v,w\})|=3d-3\mu+\nu\) and thus \(\rdeg_{G/\mathcal{P}_P}(P)=3(d-\mu)\).
	If \(\{u,v,w\}\) is not an independent set, we need to be careful to exclude the vertices themselves
	from our count. If we denote the number of edges in~\(G[\{u,v,w\}]\) by \(e\), doing this yields
	\begin{align*}
	\rdeg_{G/\mathcal{P}_P}(P)
	&=\begin{cases}
		3d-3\mu&\text{if } e=0,\\
		3d-2\mu-\lambda-2&\text{if } e=1,\\
		3d-\mu-2\lambda-3&\text{if } e=2,\\
		3d-3\lambda-3&\text{if } e=3.
	\end{cases}\\
	&\geq\frac{3}{2}\lb_1(G)-1.
	\end{align*}
	Since every contraction can decrease a red degree by at most \(1\), the number of parts of order two in the partition \(\mathcal{P}_k\) is at least \(\frac{3}{2}\lb_1(G)-1-\lb_1=\frac{\lb_1(G)}{2}-1\).
	Together with one of the parts that gets merged to form \(P\) itself, these are at least \(\nicefrac{\lb_1(G)}{2}\) disjoint contraction pairs.
	
	Now, assume that \(|P|=4\), i.e., that \(P\) is obtained by merging two pairs of vertices, say \(\{u,v\}\) and \(\{w,x\}\).
	Note that every red neighbor of \(\{u,v,w\}\) (besides possibly \(x\)) would also be a red neighbor of \(\{u,v,w,x\}\).
	Thus, the same reasoning yields that at least \(\frac{\lb_1(G)}{2}-2\) pairs disjoint from \(\{u,v,w,x\}\) must be contracted
	before \(P\) is formed. Together with the two pairs~\(\{u,v\}\) and~\(\{w,x\}\), this again yields that \(\mathcal{P}_k\)
	contains least \(\nicefrac{\lb_1}{2}\) parts of size two.
\end{proof}
\end{lemma}


In the following, we show that, extending \Cref{lem:contrac-seq-is-struct}, for many specific families of strongly regular graphs, there exists an \(\lb_1\)-contraction sequence in which (almost) all vertices are initially paired.

\subsection{Families of strongly regular graphs}
Paley and Rook's graphs are two families of strongly regular graphs known to be \(\lb_1\)-collapsible~\cite{tww_products,bounds_on_tww}. In the following, we add more families of strongly regular graphs to the list of \(\lb_1\)-collapsible graph classes, namely Johnson and Kneser graphs over 2-sets, Peisert graphs, and Latin square graphs of cyclic groups.

\subsubsection{Johnson and Kneser graphs}

For~\(n \in \N_{\geq 2}\) the Johnson graph~\(J(n,2)\) is the graph on the vertex set $\binom{[n]}{2}$ where two vertices $S$ and $S'$ in $\binom{[n]}{2}$ are adjacent precisely if~\(|S \cap S'| =  1\). Note that \(J(n,2)\) is a strongly regular graph with parameters \(\left(\binom{n}{2}, 2n-4, n-2, 4\right)\) whenever \(n \geq 4\). The two exceptions \(J(2,2)\) and \(J(3,2)\) are complete graphs on one and three vertices, respectively.
Hence
\begin{equation}\label{lem:lb1-tww-general-johnson-graphs}
	\lb_1(J(n,2))= \begin{cases}
		2(n-3) &\text{if $n \geq 5$, and}\\
		0 &\text{otherwise}.
	\end{cases} 
\end{equation}

\begin{lemma}
	\label{lem:tww-2-johnson-graphs}
	For each $n \in \mathbb{N}_{\geq 3}$ it holds that \(\tww(J(n,2)) \leq 2(n-3)\).
\end{lemma}

\begin{proof}
	Let $N \coloneqq \binom{n}{2}$ and let~\(J_r(n,2)\) be the trigraph obtained from~\(J(n,2)\) by coloring precisely the edges in~\(\{ \{u,v\}\colon n \in u\cap v\}\) and the edges incident to~\(\{n-1,n\}\) 
	 red.
	\begin{claim}\label{claim:tww-johnson-induction_step}
		For \(n \geq 3\) the graph \(J_r(n,2)\) can be contracted to an isomorphic copy of \(J_r(n-1,2)\) by a partial contraction sequence of width at most \(2(n-2)\).
	\end{claim}
	\begin{claimproof}
	Consider the following partial contraction sequence:
	In increasing order for each $i \in [n-2]$ we contract
	$\{i,n\}$ with $\{i,n-1\}$.
	Subsequently, we contract~$\{n-1,n\}$ with the vertex $\{\{n-2,n-1\}, \{n-2,n\} \}$.
	
	Let $J_r(n,2) = G_N, G_{N-1}, G_{N-2}, \dots, G_{N-(n-1)}$ be the corresponding quotient graphs.
	Fix $i \in [n-2]$.
	It holds that $\rdeg_{G_{N-i}}(\{n-1, n\}) \leq \rdeg(J_r(n,2)) = 2(n-2)$ since~$\{n-1,n\}$ is joined to all of its neighbors via red edges and all contractions that happened so far are among neighbors of~$\{n-1, n\}$.
	For~$j \in [i]$ it holds that $\rdeg_{G_{N-i}}(\{\{n-1, j\}, \{n,j\}\}) \leq 2n-5$, since there are precisely $2n-3$ vertices in $J_r(n,2)$ containing at least one of the values $n-1$ and~$n$ (observe that $\{n-1, j\}$ and $\{n,j\}$ do not account to the considered red degree). Every other vertex of $J_r(n,2)$ is adjacent to a vertex of the form~$\{k,n\}$ precisely if it is adjacent to $\{k,n-1\}$ and, hence, is of red degree~0.
	 
	In $G_{N-(n-1)}$ there are $n-2$ vertices containing $n$. These vertices form an~$(n-2)$-clique where every edge is red.
	Observe that $G_{N-(n-1)}$ is isomorphic to  $J_r(n-1,2)$.
	As argued before, the vertices $\{\{i,n-1\} \{i,n\}\}$ for $i \in [n-3]$ are not joined via red edges to vertices outside of this clique and, hence, are of red degree $n-3$.
	The vertex~$\{\{n-2,n\}, \{n-2, n-1\}, \{n-1,n\}\}$ has apart from $n-3$ red neighbors in the clique precisely the red neighbors $\{j, n-2\}$ for~$j \in [n-3]$ yielding a red degree of $2n-6$.
	All further vertices are of red degree at most~1.
	\end{claimproof}

	\begin{claim}\label{claim:tww-johnson-base_case}
		For \(n \geq 4\) the graph \(J(n,2)\) can be contracted to \(J_r(n-1,2)\) by a partial contraction sequence of width at most~\(2(n-3)\).
	\end{claim}
	\begin{claimproof}
	Since $J(n,2)$ and $J_r(n,2)$ have the same vertex set, we may
	consider the same contraction sequence as in the proof of \Cref{claim:tww-johnson-induction_step}.
	
	Let $J(n,2) = G'_{N}, G'_{N-1}, G'_{N-2}, \dots, G'_{N-(n-1)}$ be the corresponding quotient graphs of $J(n,2)$.
	For all $i \in [n-1]$ we have $V(G_{N-i}) = V(G'_{N-i})$ and
	\[\rdeg_{G'_{N-i}}(v) \leq \rdeg_{G_{N-i}}(v)~\text{for all $v \in V(G'_{N-i})$}.\] Thus, we only need to consider the red degree of those vertices and quotient graphs where the bounds in the proof of \Cref{claim:tww-johnson-induction_step} exceed $2(n-3)$.
	Fix $i \in [n-2]$.
	The red degree of $\{n-1, n\}$ in $G_{N-i}$ is $0$ since all vertices involved in merges so far are (black) neighbors of $\{n-1, n\}$ in $J(n,2)$.
	For~$j \in [i]$ it holds that~$\rdeg_{G_{N-i}}(\{\{n-1, j\}, \{n,j\}\}) \leq 2n-6$, since there are precisely~$2n-3$ vertices in $J_r(n,2)$ containing at least one of the values $n-1$ and $n$, the sets~$\{n-1,j\}$ and $\{n,j\}$ as well as $\{n-1, n\}$ do not contribute to the red degree of the considered vertex.
	\end{claimproof}
	
	The result follows by induction on $n$:
	The Johnson graph $J(3,2)$ is a complete graph on three vertices and $\tww(J(3,2)) = 0$.
	Fix $n \in \mathbb{N}_{\geq 4}$.
	By \Cref{claim:tww-johnson-base_case} there exists a width-$2(n-3)$ partial contraction sequence of~$J(n,2)$ yielding a graph isomorphic to $J_r(n-1,2)$.
	By \Cref{claim:tww-johnson-induction_step} we can contract~$J_r(n',2)$ to~$J_r(n'-1,2)$ by a partial contraction sequence of width at most~$2(n'-2)$, which is less than~$2(n-3)$ whenever $n' < n$.
	We inductively contract graphs isomorphic to~$J_r(n',2)$ to graphs isomorphic to~$J_r(n'-1,2)$ starting with~$n'= n-1$ until the resulting graph is isomorphic to $J_r(2,2)$, which is a one-vertex graph and is trivially of twin-width~0.
\end{proof}

\begin{theorem}
	\label{thm:tww-johnson-graphs}
	For each $n \in \mathbb{N}_{\geq 2}$ the Johnson graph \(J(n,2)\) satisfies
	\begin{align*}
		\tww(J(n,2))=\lb_1(J(n,2))= 
		\begin{cases}
			2(n-3) &\text{if $n \geq 5$, and}\\
			0 &\text{otherwise}.
		\end{cases}
	\end{align*}
\end{theorem}

\begin{proof}
	Since~\(J(2,2)\) and~\(J(3,2)\) are complete graphs on one and three vertices, respectively, they have twin-width 0. Further, \(J(4,2)\) is the complement of three disjoint matching edges and thus, has twin-width 0. In the case where~\(n\geq 5\), the result follows from Equation~\ref{lem:lb1-tww-general-johnson-graphs} and Lemma~\ref{lem:tww-2-johnson-graphs}.
\end{proof}

For~\(n \in \N_{\geq 2}\), analogously to Johnson graphs, the Kneser graph \(K(n,2)\) is the graph on the vertex set~\(\binom{[n]}{2}\) with two vertices being adjacent precisely if their intersection is empty.

\begin{corollary}\label{coro: Kneser}
	For each $n \in \mathbb{N}_{\geq 2}$ the Kneser graph \(K(n,2)\) satisfies
	\begin{align*}
		\tww(K(n,2))=\lb_1(K(n,2)))= 
		\begin{cases}
			2(n-3) &\text{if $n \geq 5$, and}\\
			0 &\text{otherwise}.
		\end{cases}
	\end{align*}
\end{corollary}

\begin{proof}
	This follows from Theorem~\ref{thm:tww-johnson-graphs}, since for every~\(n \in \mathbb{N}_{\geq 2}\) the Kneser graph~\(K(n,2)\) is the complement of the Johnson graph~\(J(n,2)\) and the twin-width is stable under taking complements~\cite{tww1}.
\end{proof}

\subsubsection{Self-complementary, edge- and vertex-transitive graphs}

In this subsection we generalize the result of~\cite{bounds_on_tww} that Paley graphs are \(\lb_1\)-collapsible.
We include the proof of the following folklore lemma for self-containment.

\begin{lemma}[folklore]
	\label{lem:sc+trans-conf}
	A self-complementary, edge- and vertex-transitive graph is a conference graph.
\end{lemma}

\begin{proof}
	Let $G$ be a self-complementary, edge- and vertex-transitive graph of order~$n$.
	Since $G$ is vertex-transitive and self-complementary is it $\nicefrac{(n-1)}{2}$-regular.
	From the edge-transitivity we obtain that each two adjacent vertices of $G$ have the same number of common neighbors, say $\lambda$.
	Since $G$ is self-complementary, each two non-adjacent vertices of $G$ have $n-(n-1) + \lambda =  \lambda + 1$ common neighbors.
	Altogether, $G$ is an $\srg(n, \nicefrac{(n-1)}{2}, \lambda, \lambda + 1)$.
	Applying Equation~\eqref{eq: srg} we obtain $\lambda = \nicefrac{(n-5)}{4}$ and, thus, $G$ is a conference graph.
\end{proof}

This gives a lower bound of \(\nicefrac{(n-1)}{2}\) for the twin-width of a self-complemen\-tary, edge- and vertex-transitive graph on \(n\) vertices. Next, we show that this lower bound is also already an upper bound. This generalizes the known result for the twin-width of Paley graphs and also includes the family of Peisert graphs~\cite{PEISERT01}.

\begin{theorem}\label{thm: self-complementary}
	All self-complementary, edge- and vertex-transitive graphs are \(\lb_1\)-collapsible.
\end{theorem}

\begin{proof}
	Due to a characterization of Zhang~\cite{DBLP:journals/jgt/Zhang92} all self-complementary, edge- and vertex-transitive graphs are Cayley graphs of the form~\(\Cay(\Gamma,S \cup S^{-1})\) for some abelian group~\(\Gamma\) with the property that~\(1_\Gamma \in \Gamma\) is the only self-inverse element, and some set~\(S \subset \Gamma\).
	
	We give a contraction sequence for~\(G=\Cay(\Gamma, S \cup S^{-1})\) with red degree at most~\(\frac{n-1}{2}\), where \(n=|G|\). Since~\(G\) has~\(n\) vertices, it suffices to give~\(\frac{n-1}{2}\) contraction steps without exceeding the allowed red degree of~\(\frac{n-1}{2}\).
	
	We contract all pairs \(\{x,x^{-1}\}\) for \(x \in \Gamma \setminus \{1_\Gamma\}\) in an arbitrary order and claim that this keeps the red degree bounded by \(\frac{n-1}{2}\). First of all, since \(1_\Gamma\) is the only self-inverse element, these pairs are well-defined and there are~\(\frac{n-1}{2}\) many. Further, note that after contracting \(l\in [\frac{n-1}{2}]\) pairs, an uncontracted vertex has at most red degree \(l\) (at most one red edge for each contracted pair). Hence it suffices to only consider the red degree of contracted pairs. For the first contracted pair (independently of which first contraction we choose), we have red degree~\(\frac{n-1}{2}\), since~\(G\) is a conference graph by \Cref{lem:sc+trans-conf} and thus,~\(\lb_1(G)=\frac{n-1}{2}\). From then on, the red degree of this pair can only decrease:
	Since~\(\Gamma\) is an abelian group, a vertex~\(x\) is a neighbor of~\(y\) if and only if~\(x^{-1}\) is a neighbor of~\(y^{-1}\). Thus, for every common neighbor~\(x\) of~\(y\) and~\(y^{-1}\), the inverse~\(x^{-1}\) is also a common neighbor of~\(y\) and~\(y^{-1}\).
	
	This argument also shows that for every contracted pair the red degree is less or equal to the red degree when contracting this pair first, i.e. less or equal to \(\frac{n-1}{2}\).
	Hence,~\(\tww(G) = \frac{n-1}{2}\).
\end{proof}

\subsection{Latin squares}

For each $n \in \mathbb{N}_{\geq 1}$ a \textit{Latin square} of order~\(n\) is an~\(n \times n\)-matrix over $[n]$ such that each row and each column contain each value exactly once.
To every Latin square $M$ of order~\(n\) one can associate a \textit{Latin square graph} with~\(n^2\) vertices, one vertex $(r,c)$ for each row \(r\) and column \(c\). Two distinct vertices $(r,c)$ and $(r', c')$ are adjacent precisely if either $r=r'$, $c=c'$, or $M_{r,c}=M_{r',c'}$.
The \emph{order} of a Latin square graph is the order of the underlying Latin square.

Note that a Latin square graph
depends only on the \emph{main class} of the Latin square, that is, the isomorphism type of the graph is invariant under permutations of rows, columns or symbols of the underlying Latin square,
and under switching the roles of rows, columns, or symbols.

For each $n \geq 2$ a Latin square graph \(L\) of order $n$ is an \(\srg(n^2,3(n-1),n,6) \) and, hence,
$\tww(L)\geq \lb_1(L)=\min\{4n-8,6n-9\}=4n-8$.

\subsubsection{Latin squares of cyclic groups}
Given a finite group \((G,\cdot)\), the multiplication table of \(G\) forms a Latin square. More formally, we index the rows and columns of this Latin square
by elements of \(G\) and use \(g\cdot h\) as the entry at position \((g,h)\).
We denote the corresponding Latin square graph by \(\ls(G)\).
\begin{theorem}\label{thm: latin squares}
For every \(n\in\N\), the graph \(\ls(\Z_n)\) is \(\lb_1\)-collapsible.
\begin{proof}
We first discuss that the claim is true for \(n\leq 3\).
Note that~\(\ls(\Z_1)\) and~\(\ls(\Z_2)\) are complete graphs on one and four vertices and $\ls(\Z_3)$ is the complement of the disjoint union of three 3-cliques. Hence, all three graphs are $\lb_1$-collapsible (with a lower bound of~0).
From now on we assume that~\(n\geq 4\).


\begin{figure}
\centering
\scalebox{0.96}{
\begin{tikzpicture}[yscale=-1, scale=0.5]
	\foreach \x in {0,...,10} {
		\foreach \y in {0,...,10} {
			\pgfmathtruncatemacro{\value}{mod(\x+\y,11)}
			\node at (\x,\y) {\value};
		}
	}
	\foreach \x in {1.5, 3.5, 5.5, 7.5, 9.5} {
		\draw (\x,-0.5) -- (\x,10.5);	
	}
	\foreach \y in {3.5, 7.5} {
		\draw (-0.5,\y) -- (10.5,\y);	
	}
\end{tikzpicture}
\hspace{10mm}
\begin{tikzpicture}[yscale=-1, scale=0.5]
	\foreach \x in {0,...,10} {
		\foreach \y in {0,...,10} {
			\pgfmathtruncatemacro{\value}{mod(\x+\y,11)}
			\node at (\x,\y) {\value};
		}
	}
	\foreach \x in {3.5, 5.5, 7.5, 9.5} {
		\draw (\x,-0.5) -- (\x,10.5);	
	}
	\draw (1.5,-0.5) -- (1.5, 3.5)
	      (1.5, 7.5) -- (1.5,10.5);
	\foreach \y in {3.5, 7.5} {
		\draw (-0.5,\y) -- (10.5,\y);	
	}
\end{tikzpicture}
}

\caption{The partitions \(\mathcal{P}_{1,2}\) and \(\mathcal{P}_{1,2}^{\{P_{1,0}\}}\) of the Latin square (graph) \(\ls(\Z_{11})\).}
\label{fig:latin-square_grid}
\end{figure}

Set \(\ell\coloneqq\lceil\log_2 n\rceil\).
For all \(x,y\leq\ell\), we write \(\mathcal{P}_{x,y}\) for the
partition of~\(V(\ls(\Z_n))\) which divides the Latin square into a grid of \(2^y\times 2^x\) rectangles,
that is, rectangles of height \(2^y\) and width \(2^x\), and at most one column and one row of horizontally
or vertically smaller rectangles respectively, see \Cref{fig:latin-square_grid}.

More formally, for each \(i\in\{0,\dots,\lceil \nicefrac{n}{2^y}\rceil-1\}\) and \(j\in\{0,\dots,\lceil \nicefrac{n}{2^x}\rceil-1\}\), the partition contains the part
\[P_{i,j}\coloneqq\left\{(i',j')\in\Z_n^2\colon \begin{matrix}2^yi\leq i'<\min\{2^y(i+1),n\},\\2^xj\leq j'<\min\{2^x(j+1),n\}\end{matrix}\right\},\]
where \((i',j')\) is the vertex in the \(i'\)-th row and \(j'\)-th column of the latin square graph.

A \emph{$\mathcal{P}_{x,y}$-row} is a union of all parts $P_{i,j}$ which belong to the same index~$i$. Analogously, a  \emph{$\mathcal{P}_{x,y}$-column} is a union of all parts $P_{i,j}$ which belong to the same index $j$. Two parts \(P_{i,j}\) and \(P_{i',j'}\) are \emph{horizontally consecutive} or \emph{direct horizontal neighbors}
if \(i=i'\) and \(j'=j\pm 1\bmod \lceil\nicefrac{n}{2^x}\rceil\). They are \emph{vertically consecutive} or \emph{direct vertical neighbors} if \(j=j'\)
and \(i'=i\pm 1\bmod \lceil\nicefrac{n}{2^y}\rceil\).
If \(S\subseteq\mathcal{P}_{x+1,y}\) or \(S\subseteq\mathcal{P}_{x,y+1}\), then we write $\mathcal{P}_{x,y}^S$
for the partition
\[\mathcal{P}_{x,y}^S=S\cup\{P\in \mathcal{P}_{x,y}\colon P\not\subseteq\bigcup S\}\]
obtained from \(\mathcal{P}_{x,y}\) by merging all pairs of \(\mathcal{P}_{x,y}\)-parts that together form a part in \(S\),
see \Cref{fig:latin-square_grid}.
Note that \(\mathcal{P}_{x,y} \feq \mathcal{P}_{x,y}^S\), and $\mathcal{P}_{x,y}^S \feq \mathcal{P}_{x+1,y}$ or \(\mathcal{P}_{x,y}^S \feq \mathcal{P}_{x,y+1}\). Further, if \(S\subseteq T\), then \(\mathcal{P}_{x,y}^S\feq\mathcal{P}_{x,y}^T\).

We show that every partition sequence of $\ls(\Z_n)$ containing the partitions
\[\mathcal{P}_{0,0}, \mathcal{P}_{0,1}, \mathcal{P}_{1,1},\mathcal{P}_{1,2},\dots ,\mathcal{P}_{\ell-1,\ell}, \mathcal{P}_{\ell,\ell}\]
has width at most \(4n-8\). Note that in such a partition sequence, all intermediate partitions are of the form \(\mathcal{P}_{x,x}^S\) with \(S\subseteq\mathcal{P}_{x,x+1}\) or of the form~\(\mathcal{P}_{x,x+1}^S\) with \(S\subseteq\mathcal{P}_{x+1,x+1}\). Thus, it suffices to bound the red degree in all such partitions.

\begin{claim}\label{claim:latin:P01}
For every \(S\subseteq\mathcal{P}_{0,1}\) we have \(\maxrdeg(\ls(\Z_n)/\mathcal{P}_{0,0}^S)\leq 4n-8\).
\begin{claimproof}
First, let \(P\in\mathcal{P}_{0,0}^S\) be a  singleton part.
Every red neighbor of \(P\) in $\ls(\Z_n)/\mathcal{P}_{0,0}^S$ contains a vertex that either shares a row or a symbol in $\ls(\Z_n)$ with the single vertex in \(P\).
Since there are only \(2n-2\) such vertices we obtain that \(\rdeg_{\ls(\Z_n)/\mathcal{P}_{0,0}^S}(P) \leq 2n-2\leq 4n-8\).

Now let \(P\in\mathcal{P}_{0,0}^S\) be a non-singleton part, i.e., a \(2\times 1\) rectangle.
The two vertices in \(P\) are adjacent and, hence, \(\rdeg_{\ls(\Z_n)/\mathcal{P}_{0,0}^{\{P\}}}(P)=\lb_1(\ls(\Z_n))=4n-8\).
All parts \(Q\in S\) in the same column as \(P\) are fully connected to \(P\), while all non-singleton parts \(Q\in S\) which share a row or symbol with \(P\)
contain a vertex which is joined via a red edge to \(P\) in \(\ls(\Z_n)/\mathcal{P}_{0,0}^{\{P\}}\). Altogether
\(\rdeg_{\ls(\Z_n)/\mathcal{P}_{0,0}^S}(P)\leq\rdeg_{\ls(\Z_n)/\mathcal{P}_{0,0}^{\{P\}}}(P)=4n-8\).
\end{claimproof}
\end{claim}

Claim~\ref{claim:latin:P01} already suffices to show that both \(\ls(\Z_4)\) and \(\ls(\Z_6)\) are \(\lb_1\)-collapsible:
Indeed, by the previous claim applied to even values of \(n\), we find a partial \(\lb_1\)-contraction sequence which contracts \(\ls(\Z_n)\) to at most
\(\nicefrac{n^2}{2}\) vertices. We argue that every completion of this contraction sequence
has width at most \(\lb_1\). Indeed, since the contraction sequence up to this point
has width at most \(\lb_1\), and every further contractions creates a trigraph of order at most \(\nicefrac{n^2}{2}-1\)
and thus maximum degree at most \(\nicefrac{n^2}{2}-2\), the width of the total sequence is bounded by \(\max\{4n-8,\nicefrac{n^2}{2}-2\}\).
For \(n=4\) and \(n=6\), this bound is at most \(4n-8\).
Together with the observation that \(\ls(\Z_5)\) is isomorphic to the Payley graph
on \(25\) vertices, we find that \(\ls(\Z_n)\) is \(\lb_1\)-collapsible for all \(n \leq 6\).

From now on, we thus assume that \(n\geq 7\),
which we will often implicitly use in bounds of the type \(3n-1\leq 4n-8\).

We start the main part of the proof with two special cases, namely the partitions \(\mathcal{P}_{x,y}\) themselves with \(S=\emptyset\) (\Cref{claim:latin:uncontracted}), and partitions \(\mathcal{P}_{x,y}^S\) where \(S=\{P\}\) consists of a single contracted part (\Cref{claim:latin:first_contraction}), for which
we need two technical preparations.

\begin{claim} \label{claim:latin:symbols}
If \(B\subseteq V(\ls(\Z_n))\) is an \(a\times b\) rectangle in $\ls(\Z_n)$ and $C$ is a \(\mathcal{P}_{x,y}\)-column of~$\ls(\Z_n)$, then at most $3+\frac{a+b+2^x-5}{2^y}$ \(\mathcal{P}_{x,y}\)-parts in \(C\)
share a symbol with $B$.
\end{claim}
\begin{claimproof}
Let \(r_y\in\{1,\dots,2^y\}\) be the height of the rectangles in the last row of the partition \(\mathcal{P}_{x,y}\).
Let \(s_0\) be the symbol in the top-left corner of~\(B\).
The set of symbols of \(B\) is \(S\coloneqq s_0+\{0,1,\dots,a+b-2\} \bmod n\), which is either the whole set~\(\Z_n\), or an interval
of cardinality $a+b-1$.

If \(B\) contains all symbols in \(\Z_n\), then \(a+b-1\geq n\) and thus
\[3+\frac{a+b+2^x-5}{2^y}\geq 3+\frac{n-3}{2^y}\geq\left\lceil\frac{n}{2^y}\right\rceil,\]
which yields the claim.

Otherwise, consider a \(\mathcal{P}_{x,y}\)-column $C$.
The $\mathcal{P}_{x,y}$-parts of $C$ that share a symbol with~\(B\) form an interval in the cyclic order on the \(\mathcal{P}_{x,y}\)-rows.
Assume that~$C$ consists of~\(w\) columns of~$\ls(\Z_n)$ and~\(B\) shares a symbol with~\(h\) parts of~$C$.
If \(h=1\), then we are done.
Hence, we may assume that $h>1$.
Consider the cell with symbol \(s_0\) in the last \(\ls(\Z_n)\)-column of \(C\)
and denote its row by~\(i\). Note that all cells with symbol \(s_0\) in \(C\) are contained
in rows with index in \(\{i,i+1\bmod n,\dots,i+w-1\bmod n\}\).
Similarly, the cells with symbol \(s_0+j\) are contained in the rows
with index in \(\{i+j\bmod n,i+j+1\bmod n, i+j+w-1\bmod n\}\). In total, each cell sharing a symbol with \(B\)
is contained in a row with index \(\{i,i+1\bmod n,\dots,i+a+b+w-2\bmod n\}\).

Since all \(\mathcal{P}_{x,y}\)-rows besides the exceptional row of height \(r_y\) have height~\(2^y\),
and \(B\) shares a symbol with every row of all but the outer two parts it shares a symbol with in column \(C\),
this yields that
\(2+(h-3)2^y+r_y\leq a+b+w-2\) and thus
\[h\leq\frac{a+b+w-r_y-4+3\cdot 2^y}{2^y}\leq 3+\frac{a+b+2^x-5}{2^y}.\qedhere\]
\end{claimproof}

In the following, we will repeatedly encounter expressions of the same form, hence we introduce a family of functions to describe those terms.
For two positive real numbers \(c\) and \(d\), we define the function \(f_{c,d}\colon[1,\infty)\to\mathbb{R}\) by \(f_{c,d}(z)\coloneqq \nicefrac{c}{2^z}-\nicefrac{d}{4^z}\).
\begin{claim}\label{claim:latin:functions}
If \(c\geq d\), then for every \(z\in[1,\infty)\), we have \(f_{c,d}(z)\leq \frac{c}{2}-\frac{d}{4}\).
\begin{claimproof}
Taking the derivative of \(f_{c,d}\) yields
\begin{align*}
f'_{c,d}(z)
&=-\ln(2)(2^{-z}c-2\cdot 4^{-z}d)\\
&=-\ln(2)2^{-z}(c-2^{1-z}d)\\
&\leq -\ln(2)2^{-z}(c-d)\\
&\leq 0,
\end{align*}
which proves that \(f_{c,d}\) is decreasing. Thus, \(f_{c,d}(z)\leq f_{c,d}(1)=\frac{c}{2}-\frac{d}{4}\).
\end{claimproof}
\end{claim}

First, we deal with the partitions \(\mathcal{P}_{x,y}\) themselves:
\begin{claim}\label{claim:latin:uncontracted}
For all \(x\geq 0\) and \(y\in\{x,x+1\}\), the trigraph \(\ls(\Z_n)/\mathcal{P}_{x,y}\) has maximum degree at most \(4n-8\).
\begin{claimproof}
Let \(Q\) be a part of \(\mathcal{P}_{x,y}\).
If $Q = \{q\}$ for some vertex $q \in V(\ls(\Z_n))$, then
$\deg_{\ls(\Z_n)/\mathcal{P}_{x,y}}(q) \leq  \deg_{\ls(\Z_n)}(q) = 3n-3\leq 4n-8$.
We may assume from now on that $|Q| \geq 2$ and, in particular, $y \geq 1$.

If \(x=0\) and \(y=1\), then \(Q\) is a \(2\times 1\) rectangle, which in \(\ls(\Z_n)/\mathcal{P}_{0,1}\) has exactly~\(n-1\) horizontal neighbors,
\(\lceil \nicefrac{n}{2} \rceil-1<\nicefrac{n}{2} \) vertical neighbors,
and shares a symbol with at most two parts per column and thus with at most \(2(n-1)\) parts not in its own column
in total. This counts the two direct horizontal neighbors of~\(Q\) twice. Thus, \(Q\) has a total degree of at most
\(\deg_{\mathcal{P}_{0,1}}(Q)<\frac{7}{2}n-5\leq 4n-7\). 

It remains to consider the case that \(x,y\geq 1\).
In this case \(Q\) has exactly \(\lceil\nicefrac{n}{2^y}\rceil-1<\nicefrac{n}{2^y}\) vertical and \(\lceil \nicefrac{
n}{2}^x\rceil-1<\nicefrac{n}{2^x}\) horizontal neighbors.
If \(Q\) is a rectangle of height \(1\), then \Cref{claim:latin:symbols} implies that \(Q\) shares a symbol with at most \(3+2^{x+1-y}-\nicefrac{4}{2^y}\) parts per column and thus with less than \((3+2^{x+1-y}-\nicefrac{4}{2^y})\cdot \nicefrac{n}{2^x}\) parts not in its own column.
Thus, \(Q\) has degree
\begin{equation}\label{eq:latin:height1}
\deg_{\mathcal{P}_{x,y}}(Q)<  \nicefrac{n}{2^y} + \nicefrac{n}{2^x} +  (3+2^{x+1-y}-\nicefrac{4}{2^y})\cdot \nicefrac{n}{2^x} =  (\nicefrac{4}{2^x}+\nicefrac{3}{2^y}-\nicefrac{4}{2^{(x+y)}})\cdot n.
\end{equation}
If \(x=y\), then this simplifies to \((\nicefrac{7}{2^x}-\nicefrac{4}{4^x})n=f_{7,4}(x)\cdot n\leq \frac{5}{2}n\leq 4n-7\),
where the second-to-last inequality follows from \Cref{claim:latin:functions}.
If \(x+1=y\), then \eqref{eq:latin:height1} instead simplifies to
\((\nicefrac{5.5}{2^x}-\nicefrac{2}{4^x})n=f_{\frac{11}{2},2}(x)\cdot n\leq \frac{9}{4}n\leq 4n-7\),
which again follows from \Cref{claim:latin:functions}.

If \(Q\) has height at least \(2\), then it shares a symbol with its two direct horizontal neighbors.
Moreover, \Cref{claim:latin:symbols} (with $a \leq 2^x$ and $b \leq 2^y$) implies that \(Q\) shares a symbol with at most \(4+2^{x+1-y}-\nicefrac{5}{2^y}\)  parts per \(\mathcal{P}_{x,y}\)-column.
Thus, \(Q\) shares a symbol with less than \((4+2^{x+1-y}-\nicefrac{5}{2^y})\cdot \nicefrac{n}{2^x}\) parts in \(\mathcal{P}_{x,y}\)-columns besides its own. This includes
the two horizontally directly adjacent parts of \(Q\), which we already counted. In total, we get
\begin{equation}\label{eq:latin:height_geq2}
\deg_{\mathcal{P}_{x,y}}(Q)<(\nicefrac{5}{2^x}+\nicefrac{3}{2^y}-\nicefrac{5}{2^{(x+y)}})\cdot n-2.
\end{equation}
If \(x=y\), this simplifies to \((\nicefrac{8}{2^x}-\nicefrac{5}{4^x})n-2=f_{8,5}(x)\cdot n-2\leq\frac{11}{4}n-2\leq 4n-7\).
If \(x+1=y\), then \eqref{eq:latin:height_geq2} instead simplifies to \(f_{\frac{13}{2},\frac{5}{2}}(x) \cdot n-2\leq \frac{21}{8}n-2\leq 4n-7\).
\end{claimproof}
\end{claim}

Next, we consider the red degree of a part obtained by contracting two parts of \(\mathcal{P}_{x,y}\), that is, \(\rdeg_{\ls(\Z_n)/\mathcal{P}_{x,y}^{\{P\}}}(P)\)
where \(P\in\mathcal{P}_{x+1,y}\) or \(P\in\mathcal{P}_{x,y+1}\) is a union of two parts of \(\mathcal{P}_{x,y}\).
\begin{claim}\label{claim:latin:first_contraction}
For every \(x\geq 0\) and \(y\in\{x,x+1\}\), we have
\[\rdeg_{\ls(\Z_n)/\mathcal{P}_{x,y}^{\{P\}}}(P)\leq 4n-8.\]
\begin{claimproof}
If \(x=y=0\), then the claim follows from \Cref{claim:latin:P01}.
Next, consider the partition \(\mathcal{P}_{0,1}^{\{P\}}\) obtained from \(\mathcal{P}_{0,1}\) by merging
two horizontally consecutive parts. This new part has \(n-2\) horizontal neighbors
and~\(2\cdot (\lceil \nicefrac{n}{2}\rceil-1)<n\) vertical neighbors.
If the new part has height \(1\), then it contains at most two symbols and thus
shares a symbol with at most two parts per \(\mathcal{P}_{x,y}\)-column.
In total, it thus shares a symbol with at most \(2(n-2)=2n-4\) parts in columns besides its own.
Because in this case, the horizontal neighbors of \(P\) are not red neighbors, the red degree of \(P\)
is bounded by \(3n-4\leq 4n-8\).

If the new part has height \(2\), then it contains at most three distinct symbols, say \(s-1\), \(s\) and \(s+1\).
In every column where the symbol \(s\) does not form its own part (which can only happen for a single column in the last row), \(P\) shares a symbol with at most two parts.
In the possibly exceptional single column, \(P\) shares a symbol with three parts.
In total, \(P\) thus shares a symbol with at most \(2(n-3)+3=2n-3\) parts in columns besides its own.
This count includes the two direct horizontal neighbors of our merged part, which we already counted.
Thus, the new part has red degree less than \(4n-7\).

We are left with the case that \(x\geq 1\) and \(y\in\{x,x+1\}\).
First consider~\(\mathcal{P}_{x,x}^{\{P\}}\), where \(P\) is the union of two vertically consecutive \(\mathcal{P}_{x,x}\)-parts.
Then~\(P\) has \(\lceil \nicefrac{n}{2^x}\rceil-2<\nicefrac{n}{2^x}-1\) vertical red neighbors
and \(2\cdot (\lceil \nicefrac{n}{2^x}\rceil-1)<2\cdot \nicefrac{n}{2^x}\) horizontal red neighbors.
Since~\(P\) is a rectangle of size at most \(2^{x+1}\times 2^x\), \Cref{claim:latin:symbols}
implies that~$P$ shares a symbol with at most \(5\) parts per column.
In total~$P$ shares a symbol with at most \(5\cdot \nicefrac{n}{2^x}\) parts in columns besides its own. Moreover, this counts at least six parts in the same rows but different column from \(P\). Thus
$\rdeg_{\ls(\Z_n)/\mathcal{P}_{x,x}^{\{P\}}}(P)< \nicefrac{8}{2^x}\cdot n-7\leq 4n-7$.

Finally, consider \(\mathcal{P}_{x,x+1}^{\{P\}}\) where \(P\) is the union of two horizontally consecutive parts.
As in the previous case \(P\) has less than \(\nicefrac{n}{2^x}-1\) horizontal red neighbors and less than \(2\cdot \nicefrac{n}{2^{x+1}}=\nicefrac{n}{2^x}\) vertical red neighbors.
Since \(P\) is a rectangle of size at most \(2^{x+1}\times 2^{x+1}\), \Cref{claim:latin:symbols}
implies that it shares a symbol with at most \(5\) parts per column.
In total it thus shares a symbol with at most~\(5\cdot\left(\lceil \nicefrac{n}{2^x}\rceil -2\right)<5\cdot \nicefrac{n}{2^x}-5\) parts in columns besides its own.
Moreover, this counts at least two parts in the same row but different columns from \(P\).
We obtain that $\rdeg_{\ls(\Z_n)/\mathcal{P}_{x,x+1}^{\{P\}}}(P) <\nicefrac{7}{2^x}\cdot n - 8 \leq 4n-7$.
\end{claimproof}
\end{claim}

So far, we have seen in \Cref{claim:latin:uncontracted} that the trigraphs
\(\ls(\Z_n)/\mathcal{P}_{x,y}\) with \(y\in\{x,x+1\}\) have maximum degree at most \(4n-8\),
and further in \Cref{claim:latin:first_contraction}
that for every single contracted part \(P\in\mathcal{P}_{x+1,y}\) with \(y=x+1\)
or \(P\in\mathcal{P}_{x,y+1}\) with \(y=x\), we have
\(\rdeg_{\ls(\Z_n)/\mathcal{P}_{x,y}^{\{P\}}}(P)\leq 4n-8\).
Finally, we will bound the red degree in the general partitions \(\ls(\Z_n)/\mathcal{P}_{x,y}^S\), which then yields the claim of the theorem.

Let \(P\in \mathcal{P}_{x,y}^{S}\).
If \(P\in\mathcal{P}_{x,y}\), that is, \(P\) was not contracted from \(\mathcal{P}_{x,y}\) to \(\mathcal{P}_{x,y}^S\),
then
\[\rdeg_{\ls(\Z_n)/\mathcal{P}_{x,y}^S}(P)\leq \deg_{\ls(\Z_n)/\mathcal{P}_{x,y}^S}(P)\leq \deg_{\ls(\Z_n)/\mathcal{P}_{x,y}}(P),\]
since every contraction of two parts disjoint from \(P\) does not increase the degree of \(P\).
Thus, the bound \(\rdeg_{\ls(\Z_n)/\mathcal{P}_{x,y}^S}(P)\leq 4n-8\) follows from \Cref{claim:latin:uncontracted}.

If \(P\notin \mathcal{P}_{x,y}\), that is, if \(P\) is the contraction of two parts in \(\mathcal{P}_{x,y}\), then
\[\rdeg_{\ls(\Z_n)/\mathcal{P}_{x,y}^S}(P)\leq\rdeg_{\ls(\Z_n)/\mathcal{P}_{x,y}^{\{P\}}}(P).\]
Indeed, a contraction of two parts \(Q,Q'\in\mathcal{P}_{x,y}\), at least one of which was already connected via a red edge to \(P\)  does not increase the red degree of \(P\).
If both parts are not joined to \(P\) via red edges but their contraction \(Q\cup Q'\) is,
then one of \(Q\) or \(Q'\) is fully connected to~\(P\) while the other part and \(P\) are disconnected.
Since~\(|P|>1\), it contains multiple different symbols which forbids \(P\) from being fully connected to any part it
does not share a row or column with. But if \(Q\) or \(Q'\) shares a row or column with \(P\), then both parts share a row or column with \(P\),
since otherwise, they would not be contracted to reach the next partition \(\mathcal{P}_{x+1,y}\) or \(\mathcal{P}_{x,y+1}\).
Hence, both \(Q\) and \(Q'\) share an edge with \(P\), and their contraction does
not increase the red degree of \(P\).
Thus, the bound \(\rdeg_{\ls(\Z_n)/\mathcal{P}_{x,y}^S}(P)\leq 4n-8\) follows from \Cref{claim:latin:first_contraction}.

This concludes the proof that \(\maxrdeg(\ls(\Z_n)/\mathcal{P}_{x,y}^S)\leq 4n-8\), which proves that \(\ls(\Z_n)\) is \(\lb_1\)-collapsible.
\end{proof}
\end{theorem}
A similar technique of contracting vertices into rectangular parts also works to show that \(\ls(\Z_2^k)\) is \(\lb_1\)-collapsible for every \(k\),
by contracting~\(\ls(\Z_2^k)\) to~\(\ls(\Z_2^{k-1})\) inductively. More generally, by contracting the cyclic subsquares in a larger latin square \(\ls(A)\)
for an abelian group, one can show that latin square graphs of all abelian groups are \(\lb_1\)-collapsible.

\subsubsection{Upper bounds for general Latin squares}
We obtain an upper bound on the twin-width of a Latin square graph~\(G\) obtained from an \(n \times n\) Latin square from Theorem~\ref{thm:degeneracy},
which implies that \(\tww(G)\leq \sqrt{6}n^{\nicefrac{3}{2}}+O(n)\). We deem it plausible that this upper bound of \(\Theta(n^{\nicefrac{3}{2}})\) is asymptotically tight.
In particular, it would follow from a conjecture of Linial and Luria \cite{LinialLuria} that there exist Latin square graphs such that every subset of~\(\Theta(\sqrt{n})\) many vertices has a neighborhood of size~\(\Theta(n^{\nicefrac{3}{2}})\). This greatly restricts the possible partitions of these Latin squares graphs having small red degree.
\begin{conjecture}\label{conj:tww_of_random_ls}
	 The twin-width of a random Latin square graph $G$ is asymptotically larger than $\lb_1(G)$.
\end{conjecture}

\subsection{Computational results}
\label{sec: computational results}

We end with computational results on the twin-width of strongly regular graphs.
\begin{lemma}\label{lem:srg:computations}
The following classes of strongly regular graphs are \(\lb_1\)-collapsible:
\begin{itemize}
\item strongly regular graphs of order at most \(36\),
\item Latin square graphs of order at most \(8\) (i.e., with at most $64$ vertices),
\item intercalate-free Latin square graphs of order at most \(9\) (i.e., with at most~81 vertices),
\item all strongly regular graphs given by Spence in \cite{HPSpence} as of March 2025,
\item the 31\,490\,375 strongly regular graphs with parameters $(57, 24, 11, 9)$ given by Ihringer in \cite[\S3.2]{Ihringer2022}
\item the 13\,505\,292 strongly regular graphs with parameters $(63, 30, 13, 15)$ given by Ihringer in \cite[\S3.3]{Ihringer2022},
\item the 16\,565\,438 strongly regular graphs with parameters $(81, 30, 9, 12)$ given by Ihringer in \cite[\S3.9]{Ihringer2022}.
\end{itemize}
\begin{proof}
All of these classes are completely enumerated:
Strongly regular graphs of order at most \(36\) were completely enumerated by a variety of authors, see~\cite{list_of_srgs}.
The corresponding graph data is available as part of the classification
of association schemes of order up to \(34\)~\cite{list_of_as}, where
the missing thin schemes only contain the strongly regular graphs \(nK_2\)
and \(\overline{nK_2}\), which are \(lb_1\)-collapsible. The strongly regular
graphs of order \(35\) and \(36\) were classified by Spence and McKay and
are available at \cite{HPSpence}.

The (main classes of) Latin squares of order at most \(8\) and intercalate-free Latin squares of order at most \(9\) are enumerated and made available by McKay~\cite{list_of_ls}.

The strongly regular graphs of order at least \(37\) listed by Spence in~\cite{HPSpence} contains as of March 2025:
\begin{itemize}
\item a partial list of 6760 strongly regular graphs with parameters \((37,18,8,9)\),
\item the \(28\) strongly regular graphs with parameters \((40,12,2,4)\)~\cite{spence_srg_40_12_2_4},
\item the \(78\) strongly regular graphs with parameters \((45,12,3,3)\)~\cite{spence_srg_45_12_3_3},
\item the unique strongly regular graph with parameters \((49,12,5,2)\), the \(7\times 7\) Rook's graph,
\item a partial list of \(18\) strongly regular graphs with parameters \((50,21,8,9)\),
\item the \(167\) strongly regular graphs with parameters \((64,18,2,6)\)~\cite{spence_srg_64_18_2_6}.
\end{itemize}

The classes of graph described by Ihringer in \cite{Ihringer2022} are available
at~\cite{HPIhringer}.

The twin-width of all these graphs was computed by version 0.0.3-SNAPSHOT of the hydraprime solver~\cite{hydraprime}.
\end{proof}
\end{lemma}

		\section{Conclusion and further research}
We showed that extremal graphs of bounded degree and high twin-width are asymmetric (\Cref{cor:extremal_trigraphs_asymmetric}, \Cref{thm: extremal trigraphs of bd degree}). Further, we proved a new upper bound on the twin-width of degenerate graphs (Theorem~\ref{thm:degeneracy}). This bound is tight whenever the degeneracy of the considered graph is in $\omega(\log n)$, where $n$ is the order of the graph. It remains open whether
the bound is also tight in the case of graphs with small degeneracy.

Moreover, we proved the following classes of strongly regular graphs to be $\lb_1$-collapsible: Johnson graphs over 2-sets (Theorem~\ref{thm:tww-johnson-graphs}), Kneser graphs over 2-sets (Corollary~\ref{coro: Kneser}), self-complementary vertex- and edge-transitive graphs (Theorem~\ref{thm: self-complementary}), and Latin square graphs of cyclic groups (Theorem~\ref{thm: latin squares}). 
For all of these graph classes, our proofs are constructive, that is, we provide $\lb_1$-contraction sequences.
By computer search (see Section~\ref{sec: computational results}), we found more than~60~million strongly regular graphs which are $\lb_1$-collapsible and not even one such graph which is not $\lb_1$-collapsible.
However, due to the pseudo-random behavior of large strongly regular graphs,  
we believe the following, which generalizes \Cref{conj:tww_of_random_ls}:
\begin{conjecture}
	Almost all strongly regular graph are not $\lb_1$-collapsible.
\end{conjecture}
Therefore, an important future challenge is to first find a single strongly regular graph that is not $\lb_1$-collapsible.
In contrast to this, a large class of strongly regular graphs which are candidates for being $\lb_1$-collapsible are \mbox{rank-3} graphs (which are graphs adhering to a very strong symmetry condition, see also~\cite{BambergGLR23} and Chapter~6 of~\cite{cameron2004strongly}).

Finally, considering the bounds on the \(\lb_1\) in \Cref{thm:bounds_on_lb1},
and the fact that all symmetric self-complementary graphs are \(\lb_1\)-collapsible,
we conjecture the following, parts of which was already posed as a question in~\cite{bounds_on_tww}:
\begin{conjecture}
	Every simple graph of order $n$ satisfies $\tww(G) \leq \frac{n-1}{2}$ with equality if and only if $G$ is a conference graph.
\end{conjecture}

\bigskip

\noindent
\textbf{Acknowledgements.}
The authors would like to thank the anonymous reviewers, whose feedback
significantly improved the quality of the paper.
The first author’s research leading to these results has received funding
from the European Research Council (ERC) under the European Union’s
Horizon 2020 research and innovation programme (EngageS: grant agreement
No. 820148).

		\bibliography{bibliography.bib}

\begin{thebibliography}{10}

\bibitem{tww_random_graphs}
Jungho Ahn, Debsoumya Chakraborti, Kevin Hendrey, Donggyu Kim, and Sang-il Oum.
\newblock Twin-width of random graphs.
\newblock {\em Random Structures \& Algorithms}, 65(4):794--831, 2024.
\newblock \href {https://doi.org/10.1002/rsa.21247}
  {\path{doi:10.1002/rsa.21247}}.

\bibitem{bounds_on_tww}
Jungho Ahn, Kevin Hendrey, Donggyu Kim, and Sang-il Oum.
\newblock Bounds for the {{Twin-Width}} of {{Graphs}}.
\newblock {\em SIAM Journal on Discrete Mathematics}, 36(3):2352--2366, 2022.
\newblock \href {https://doi.org/10.1137/21M1452834}
  {\path{doi:10.1137/21M1452834}}.

\bibitem{finite_group_theory}
Michael Aschbacher.
\newblock {\em Finite group theory}, volume~10 of {\em Cambridge Studies in
  Advanced Mathematics}.
\newblock Cambridge University Press, Cambridge, second edition, 2000.
\newblock \href {https://doi.org/10.1017/CBO9781139175319}
  {\path{doi:10.1017/CBO9781139175319}}.

\bibitem{BambergGLR23}
John Bamberg, Michael Giudici, Jesse Lansdown, and Gordon~F. Royle.
\newblock Separating rank 3 graphs.
\newblock {\em Eur. J. Comb.}, 112:103732, 2023.
\newblock \href {https://doi.org/10.1016/J.EJC.2023.103732}
  {\path{doi:10.1016/J.EJC.2023.103732}}.

\bibitem{bollobas_maxdegree}
B{\'{e}}la Bollob{\'{a}}s.
\newblock Degree sequences of random graphs.
\newblock {\em Discret. Math.}, 33(1):1--19, 1981.
\newblock \href {https://doi.org/10.1016/0012-365X(81)90253-3}
  {\path{doi:10.1016/0012-365X(81)90253-3}}.

\bibitem{tww8_conference}
{\'{E}}douard Bonnet, Dibyayan Chakraborty, Eun~Jung Kim, Noleen K{\"{o}}hler,
  Raul Lopes, and St{\'{e}}phan Thomass{\'{e}}.
\newblock Twin-width {VIII:} delineation and win-wins.
\newblock In Holger Dell and Jesper Nederlof, editors, {\em 17th International
  Symposium on Parameterized and Exact Computation, {IPEC} 2022, September 7-9,
  2022, Potsdam, Germany}, volume 249 of {\em LIPIcs}, pages 9:1--9:18. Schloss
  Dagstuhl - Leibniz-Zentrum f{\"{u}}r Informatik, 2022.
\newblock \href {https://doi.org/10.4230/LIPICS.IPEC.2022.9}
  {\path{doi:10.4230/LIPICS.IPEC.2022.9}}.

\bibitem{tww3_conference}
{\'{E}}douard Bonnet, Colin Geniet, Eun~Jung Kim, St{\'{e}}phan Thomass{\'{e}},
  and R{\'{e}}mi Watrigant.
\newblock Twin-width {III:} max independent set, min dominating set, and
  coloring.
\newblock In Nikhil Bansal, Emanuela Merelli, and James Worrell, editors, {\em
  48th International Colloquium on Automata, Languages, and Programming,
  {ICALP} 2021, July 12-16, 2021, Glasgow, Scotland (Virtual Conference)},
  volume 198 of {\em LIPIcs}, pages 35:1--35:20. Schloss Dagstuhl -
  Leibniz-Zentrum f{\"{u}}r Informatik, 2021.
\newblock \href {https://doi.org/10.4230/LIPIcs.ICALP.2021.35}
  {\path{doi:10.4230/LIPIcs.ICALP.2021.35}}.

\bibitem{tww2}
{\'{E}}douard Bonnet, Colin Geniet, Eun~Jung Kim, St{\'{e}}phan Thomass{\'{e}},
  and R{\'{e}}mi Watrigant.
\newblock Twin-width {{II}}: Small classes.
\newblock {\em Combinatorial Theory}, 2(2), 2022.
\newblock \href {https://doi.org/10.5070/C62257876}
  {\path{doi:10.5070/C62257876}}.

\bibitem{tww7_arXiv}
{\'{E}}douard Bonnet, Colin Geniet, Romain Tessera, and St{\'{e}}phan
  Thomass{\'{e}}.
\newblock Twin-width {{VII}}: Groups.
\newblock 2022.
\newblock \href {https://arxiv.org/abs/2204.12330} {\path{arXiv:2204.12330}}.

\bibitem{tww4}
{\'{E}}douard Bonnet, Ugo Giocanti, Patrice Ossona~de Mendez, Pierre Simon,
  St{\'{e}}phan Thomass{\'{e}}, and Szymon Torunczyk.
\newblock Twin-width {IV:} ordered graphs and matrices.
\newblock {\em J. {ACM}}, 71(3):21, 2024.
\newblock \href {https://doi.org/10.1145/3651151} {\path{doi:10.1145/3651151}}.

\bibitem{tww1}
{\'{E}}douard Bonnet, Eun~Jung Kim, St{\'{e}}phan Thomass{\'{e}}, and
  R{\'{e}}mi Watrigant.
\newblock Twin-width {I:} tractable {FO} model checking.
\newblock {\em J. {ACM}}, 69(1):3:1--3:46, 2022.
\newblock \href {https://doi.org/10.1145/3486655} {\path{doi:10.1145/3486655}}.

\bibitem{minibaum}
Gunnar Brinkmann.
\newblock Fast generation of cubic graphs.
\newblock {\em J. Graph Theory}, 23(2):139--149, 1996.
\newblock \href
  {https://doi.org/10.1002/(SICI)1097-0118(199610)23:2\%3C139::AID-JGT5\%3E3.0.CO;2-U}
  {\path{doi:10.1002/(SICI)1097-0118(199610)23:2\%3C139::AID-JGT5\%3E3.0.CO;2-U}}.

\bibitem{snarkhunter}
Gunnar Brinkmann, Jan Goedgebeur, and Brendan~D. McKay.
\newblock Snarkhunter.
\newblock URL: \url{https://caagt.ugent.be/cubic/}.

\bibitem{list_of_srgs}
Andries~E. Brouwer.
\newblock Parameters of strongly regular graphs.
\newblock Accessed at: 2025-03.
\newblock URL: \url{https://aeb.win.tue.nl/graphs/srg/srgtab.html}.

\bibitem{srgBrouwer}
Andries~E. Brouwer and Hendrik Van~Maldeghem.
\newblock {\em {Strongly regular graphs}}, volume {182}.
\newblock {Cambridge University Press}, {2022}.

\bibitem{cameron2004strongly}
Peter~J Cameron.
\newblock Strongly regular graphs.
\newblock {\em Topics in algebraic graph theory}, 102:203--221, 2005.

\bibitem{spence_srg_45_12_3_3}
Kris Coolsaet, Jan Degraer, and Edward Spence.
\newblock The strongly regular (45, 12, 3, 3) graphs.
\newblock {\em Electron. J. Comb.}, 13(1), 2006.
\newblock \href {https://doi.org/10.37236/1058} {\path{doi:10.37236/1058}}.

\bibitem{house_of_graphs_cubic}
Kris Coolsaet, Sven D’hondt, and Jan Goedgebeur.
\newblock Connected cubic graphs.
\newblock {Accessed at: 03-2025}.
\newblock URL: \url{https://houseofgraphs.org/meta-directory/cubic}.

\bibitem{house_of_graphs}
Kris Coolsaet, Sven D’hondt, and Jan Goedgebeur.
\newblock House of graphs 2.0: A database of interesting graphs and more.
\newblock {\em Discrete Applied Mathematics}, 325:97--107, 2023.
\newblock {Available at \url{https://houseofgraphs.org}}.
\newblock \href {https://doi.org/10.1016/j.dam.2022.10.013}
  {\path{doi:10.1016/j.dam.2022.10.013}}.

\bibitem{guthm_software}
Holger Dell, Anselm Haak, Frank Kammer, Alexander Leonhardt, Johannes Meintrup,
  Meyer Ulrich, and Manuel Penschuck.
\newblock Guthm and guthmi: Exact and heurstic twin-width solvers, June 2023.
\newblock \href {https://doi.org/10.5281/zenodo.7996074}
  {\path{doi:10.5281/zenodo.7996074}}.

\bibitem{DBLP:books/daglib/0037866}
Christopher~D. Godsil and Gordon~F. Royle.
\newblock {\em Algebraic Graph Theory}.
\newblock Graduate texts in mathematics. Springer, 2001.

\bibitem{spence_srg_64_18_2_6}
Willem~H. Haemers and Edward Spence.
\newblock The pseudo-geometric graphs for generalized quadrangles of order (3,
  t).
\newblock {\em European Journal of Combinatorics}, 22(6):839--845, 2001.
\newblock \href {https://doi.org/10.1006/eujc.2001.0507}
  {\path{doi:10.1006/eujc.2001.0507}}.

\bibitem{list_of_as}
Akihide Hanaki.
\newblock Classification of small association schemes, January 2020.
\newblock \href {https://doi.org/10.5281/zenodo.3627821}
  {\path{doi:10.5281/zenodo.3627821}}.

\bibitem{tww_sparse_random_graphs}
Kevin Hendrey, Sergey Norin, Raphael Steiner, and Jérémie Turcotte.
\newblock Twin-width of sparse random graphs.
\newblock {\em Combinatorics, Probability and Computing}, pages 1--20, 2024.
\newblock \href {https://doi.org/10.1017/S0963548324000439}
  {\path{doi:10.1017/S0963548324000439}}.

\bibitem{HPIhringer}
Ferdinand Ihringer.
\newblock Many srgs.
\newblock Accessed at: 2025-03.
\newblock URL: \url{https://math.ihringer.org/srgs.php}.

\bibitem{Ihringer2022}
Ferdinand Ihringer.
\newblock Switching for small strongly regular graphs.
\newblock {\em Australasian Journal of Combinatorics}, 84(1):28--48, 2022.

\bibitem{guthm}
Alexander Leonhardt, Holger Dell, Anselm Haak, Frank Kammer, Johannes Meintrup,
  Ulrich Meyer, and Manuel Penschuck.
\newblock Pace solver description: Exact (guthmi) and heuristic (guthm).
\newblock In Neeldhara Misra and Magnus Wahlstr\"{o}m, editors, {\em 18th
  International Symposium on Parameterized and Exact Computation (IPEC 2023)},
  volume 285 of {\em Leibniz International Proceedings in Informatics
  (LIPIcs)}, pages 37:1--37:7, Dagstuhl, Germany, 2023. Schloss Dagstuhl --
  Leibniz-Zentrum f{\"u}r Informatik.
\newblock \href {https://doi.org/10.4230/LIPIcs.IPEC.2023.37}
  {\path{doi:10.4230/LIPIcs.IPEC.2023.37}}.

\bibitem{LinialLuria}
Nathan Linial and Zur Luria.
\newblock Discrepancy of high-dimensional permutations.
\newblock {\em Discrete Anal.}, pages Paper No. 11, 8, 2016.
\newblock \href {https://doi.org/10.19086/da.845} {\path{doi:10.19086/da.845}}.

\bibitem{list_of_ls}
Brendan~D. McKay.
\newblock Latin squares.
\newblock Accessed at: 2025-03.
\newblock URL: \url{https://users.cecs.anu.edu.au/~bdm/data/latin.html}.

\bibitem{HPnauty}
Brendan~D. McKay and Adolfo Piperno.
\newblock nauty and {T}races.
\newblock Accessed at: 2025-03.
\newblock URL: \url{https://pallini.di.uniroma1.it/}.

\bibitem{nauty}
Brendan~D. McKay and Adolfo Piperno.
\newblock Practical graph isomorphism, {II}.
\newblock {\em J. Symb. Comput.}, 60:94--112, 2014.
\newblock \href {https://doi.org/10.1016/J.JSC.2013.09.003}
  {\path{doi:10.1016/J.JSC.2013.09.003}}.

\bibitem{genreg}
Markus Meringer.
\newblock Fast generation of regular graphs and construction of cages.
\newblock {\em J. Graph Theory}, 30(2):137--146, 1999.
\newblock \href
  {https://doi.org/10.1002/(SICI)1097-0118(199902)30:2\%3C137::AID-JGT7\%3E3.0.CO;2-G}
  {\path{doi:10.1002/(SICI)1097-0118(199902)30:2\%3C137::AID-JGT7\%3E3.0.CO;2-G}}.

\bibitem{mertens}
Franz Mertens.
\newblock Ein {B}eitrag zur analytischen {Z}ahlentheorie.
\newblock {\em J. Reine Angew. Math.}, 78:46--62, 1874.
\newblock \href {https://doi.org/10.1515/crll.1874.78.46}
  {\path{doi:10.1515/crll.1874.78.46}}.

\bibitem{hydraprime}
Yosuke Mizutani, David Dursteler, and Blair~D. Sullivan.
\newblock Pace solver description: Hydra prime.
\newblock In Neeldhara Misra and Magnus Wahlstr\"{o}m, editors, {\em 18th
  International Symposium on Parameterized and Exact Computation (IPEC 2023)},
  volume 285 of {\em Leibniz International Proceedings in Informatics
  (LIPIcs)}, pages 36:1--36:5, Dagstuhl, Germany, 2023. Schloss Dagstuhl --
  Leibniz-Zentrum f{\"u}r Informatik.
\newblock \href {https://doi.org/10.4230/LIPIcs.IPEC.2023.36}
  {\path{doi:10.4230/LIPIcs.IPEC.2023.36}}.

\bibitem{PEISERT01}
Wojciech Peisert.
\newblock All self-complementary symmetric graphs.
\newblock {\em Journal of Algebra}, 240(1):209--229, 2001.
\newblock \href {https://doi.org/https://doi.org/10.1006/jabr.2000.8714}
  {\path{doi:https://doi.org/10.1006/jabr.2000.8714}}.

\bibitem{tww_products}
William Pettersson and John Sylvester.
\newblock Bounds on the twin-width of product graphs.
\newblock {\em Discret. Math. Theor. Comput. Sci.}, 25(1), 2023.
\newblock \href {https://doi.org/10.46298/DMTCS.10091}
  {\path{doi:10.46298/DMTCS.10091}}.

\bibitem{pyber_pseudorandom}
László Pyber.
\newblock {Large connected strongly regular graphs are Hamiltonian}.
\newblock 2014.
\newblock \href {https://arxiv.org/abs/1409.3041} {\path{arXiv:1409.3041}}.

\bibitem{SAT_approach}
Andr{\'{e}} Schidler and Stefan Szeider.
\newblock A {SAT} approach to twin-width.
\newblock In Cynthia~A. Phillips and Bettina Speckmann, editors, {\em
  Proceedings of the Symposium on Algorithm Engineering and Experiments,
  {ALENEX} 2022, Alexandria, VA, USA, January 9-10, 2022}, pages 67--77.
  {SIAM}, 2022.
\newblock \href {https://doi.org/10.1137/1.9781611977042.6}
  {\path{doi:10.1137/1.9781611977042.6}}.

\bibitem{spence_srg_40_12_2_4}
Edward Spence.
\newblock The strongly regular (40, 12, 2, 4) graphs.
\newblock {\em Electron. J. Comb.}, 7, 2000.
\newblock \href {https://doi.org/10.37236/1500} {\path{doi:10.37236/1500}}.

\bibitem{HPSpence}
Edward Spence.
\newblock Strongly regular graphs on at most 64 vertices, 2025.
\newblock Accessed at: 2025-03.
\newblock URL: \url{https://www.maths.gla.ac.uk/~es/srgraphs.php}.

\bibitem{DBLP:journals/jgt/Zhang92}
Hong Zhang.
\newblock Self-complementary symmetric graphs.
\newblock {\em J. Graph Theory}, 16(1):1--5, 1992.

\end{thebibliography}
		\bibliographystyle{plainurl}
\end{document}